\documentstyle[amscd,amssymb,epsf,10pt]{amsart}
\newtheorem{thm}{Theorem}[section]
\newtheorem{lem}[thm]{Lemma}
\newtheorem{conj}[thm]{Conjecture}
\newtheorem{prop}[thm]{Proposition}
\newtheorem{cor}[thm]{Corollary}

\theoremstyle{definition}
\newtheorem{defn}[thm]{Definition}

\newtheorem{rem}[thm]{Remark}

\makeatletter
\def\theequation{\@arabic\c@equation}
\makeatother
\numberwithin{equation}{thm}

\def\ff{{\Bbb F}}
\def\cc{{\cal C}}

\def\SS{{\cal S}}
\def\eps{\varepsilon}
\def\R{{\Bbb R}}

\def\G{{\cal G}}

\def\Q{{\Bbb Q}}
\def\Z{{\Bbb Z}}

\def\P{{\Bbb P}}
\def\C{{\Bbb C}}

\def\xx{{\Bbb X}}
\def\Quot{{\hbox{\rom{Quot}}}}

\def\E{{\cal E}}

\def\oo{\OO}

\def\eps{\varepsilon}

\def\whsq{\vbox to 5.8pt
{\offinterlineskip\hrule
\hbox to 5.8pt{\vrule height
5.1pt\hss\vrule height 5.1pt}\hrule}}

\def\maps{\longrightarrow}
\def\oo{{\cal O}}

\def\H{{\frak H}}
\def\C{{\Bbb C}}
\def\R{{\Bbb R}}
\def\P{{\Bbb P}}

\def\Z{{\Bbb Z}}
\def\Q{{\Bbb Q}}

\def\eps{\varepsilon}

\def\phi{\varphi}

\def\({\left(}
\def\){\right)}

\def\part{P(n)}

\def\<{\langle}
\def\>{\rangle}

\def\ee{{\Bbb E}}

\begin{document}

\title[Theta functions and Hodge numbers of moduli spaces]
{Theta functions and Hodge numbers of moduli spaces of  sheaves on rational
surfaces}

\author{Lothar G\"ottsche}
\address{
International Center for Theoretical Physics, Strada Costiera 11, P.O.
Box 586, \\ 34100 Trieste, Italy}
\email{gottsche@@ictp.trieste.it}
\maketitle\


\tableofcontents

\section{Introduction}

Let $(S,H)$ be a rational algebraic surface with an ample divisor. We assume that
$K_SH\le 0$.  In the current paper we want to compute the Betti numbers and Hodge
numbers  of the moduli spaces $M_S^H(C,d)$ of $H$-semistable 
 torsion-free sheaves of rank $2$ on $S$. 

In \cite{V-W} Vafa and Witten made a number of predictions about the
Euler numbers of moduli spaces of  sheaves on algebraic surfaces:
in many cases their generating functions should be given by modular forms.
In the case of rational surfaces this cannot be true for all 
polarizations $H$: The moduli spaces and their Euler numbers depend on 
$H$, and this dependence is not compatible with the modularity properties. 
We study the limit of the generating function for the Euler numbers as $H$ 
approaches a point $F$ on the boundary of the ample cone with $F^2=0$
(see below for the definitions).
It turns out that this limit is indeed a (quasi)-modular form
(see section 2.3).

More generally we will relate the generating functions for 
the Hodge numbers
and Betti numbers of the $M_S^H(C,d)$ to certain theta functions 
of indefinite lattices, which were introduced and studied in 
\cite{G-Z} in order to show structural results about Donaldson invariants.
That the Euler numbers and signatures are given by modular
and quasimodular forms follows then from the fact that these 
theta functions are Jacobi forms.
As in \cite{G-Z}, where the Donaldson invariants were studied,
the theta functions enter the calculations by summing over
walls.  The ample cone has a chamber structure, and the 
moduli spaces $M^H_S(C,d)$ only change  when $H$ crosses a wall.
The structure of the walls for the moduli spaces is precisely
the same as for the Donaldson invariants. Therefore we can use again
the same theta functions as in \cite{G-Z}.
We write our results for the $\chi_y$-genera instead of for the Hodge numbers,
which is equivalent as all the cohomology is of type $(p,p)$ \cite{Be}.
One could also have instead used the Poincar\'e polynomial, but I believe
that in general the $\chi_y$-genus will be better behaved.
By specializing the generating functions for the  $\chi_y$-genera of the 
moduli spaces, we also 
obtain  that the generating functions for  the signatures are given by 
modular forms, a
fact that does not seem to  have been predicted by the physics literature. 
It turns out
that the generating function for the  signatures is  better behaved than 
that for the
Euler numbers. If $F$ lies on the boundary of the positive cone, then the 
corresponding 
generating function for the signatures is a modular form and not just a 
quasimodular
form.

A surprising and interesting result is that the signatures of the 
moduli spaces $M^H_S(C,d)$ are closely related to the 
corresponding Donaldson invariants $\Phi_C^{S,H}$.
For any point $H$ in the ample cone, the generating function for the 
signatures is also the generating function for the 
Donaldson invariants $\Phi^{S,H}_C(p^r)$ evaluated on all powers of
the point class $p\in H_0(S,\Z)$. 
The signatures of the moduli spaces are just the 
coefficients of the Fourier development of this generating function,
whereas the Donaldson invariants are (up to some elementary factors)
the coefficients of the development of this function 
 into powers of a
modular function $\overline u(\tau)$ for $\Gamma(2)$.
In particular knowing all the signatures of the moduli spaces 
  $M^H_S(C,d)$ is equivalent to knowing all the 
Donaldson invariants $\Phi^{S,H}_C(p^r)$.
This relation also persists under our extension of the generating functions
and, together with the formulas for the $K3$ surfaces, suggests a similar result
for any algebraic surface.
The proof of this result uses the conjecture of Kotschick and Morgan
\cite{K-M}. Feehan and Leness \cite{F-L1}, \cite{F-L2}, \cite{F-L3}, 
\cite{F-L4} are working towards the proof of this conjecture.

This paper grew out of discussions with
Jun Li on some aspects of \cite{V-W}.
I would like to thank K.~Yoshioka for several very useful comments,
 G.~Thompson for useful discussions
and  the referee for many useful comments and improvements.

While preparing this manuscript I learned about related work.
In \cite{M-N-V-W} new  predictions are made about the
Euler numbers of $M^F_S(C,d)$, where $S$ is an rational elliptic surface, $F$ is the
class of a fibre and $CF$ even.  Yoshioka
\cite{Y4} has shown these predictions. Li
and Qin (\cite{L-Q1}, \cite{L-Q2}) have shown  blowup formulas for the 
Euler numbers and
virtual Hodge polynomials of
$M_S^H(C,d)$ for arbitrary $S$. 
After this paper was submitted Baranovsky \cite{Ba} displayed an action of
the oscilator algebra on the cohomology of the moduli spaces $M^F_S(r,C,d)$
and gave a simple relation between the Betti numbers of the
Gieseker and  Uhlenbeck compactifications. 

\section{Notations, definitions and background}

In this paper 
$S$ usually denotes  a smooth  algebraic surface over $\C$.
Often we will assume $S$ to be also rational.
For a variety $Y$ over $\C$, we denote by upper case letters the classes in $H^2(Y,\C)$,
unless they appear as walls (see below), when we denote them by Greek letters.  For
$A,B\in H^2(Y,\C)$  the intersection product on $H^2(Y,\C)$  is just denoted by   $AB$.
Later we will also need the negative of the intersection product, which we 
denote by $\<A,B\>$.
For a smooth compact variety $Y$ of complex dimension $d$ let
$$h(Y,x,y):=\sum_{p,q}(-1)^{p+q}h^{p,q}(Y)x^py^q$$ be the Hodge polynomial
(note   the   signs), and 
let 
$$H(Y)=H(Y:x,y):=(xy)^{-\frac{d}{2}}h(Y,x,y).$$  
The advantage of this (Laurent) polynomial in $x^{1/2},y^{1/2}$ is that it is
symmetric around degree $0$. In a similar way let    
$P(Y)=P(Y:y)=\sum_i (-1)^ib^i(Y) y^{i-d}:=H(Y:y,y)$ be the
(shifted) Poincar\'e polynomial (again note the signs) and let
$X_{y}(Y)=H(Y:1,y)$ be the (shifted) $\chi_{-y}$-genus.  
Then the Euler number of $Y$ is $e(Y)=X_1(Y)=P(Y,1)$, and the signature is  
$\sigma(Y):=(-1)^{\frac{d}{2}}X_{-1}(Y)$.

\subsection{Virtual Hodge polynomials and the Weil conjectures}

Virtual Hodge polynomials were
introduced in \cite{D-K}.
For  $Y$ a complex variety  
the cohomology $H^k_c(Y,\Q)$ with compact support carries
 a natural mixed 
Hodge structure. If $Y$ is smooth and projective,
this Hodge structure coincides with the classical one.
Following \cite{Ch},  we put
\begin{align*} h_v(Y:x,y)&:= 
\sum_{p,q} \sum_{k} (-1)^k h^{p,q}(H^k_c(Y,\Q)) x^py^q.
\end{align*}

These virtual Hodge polynomials have the following properties (see \cite{Ch}).
If $Y$ is a smooth projective variety, then $h_v(Y:x,y)=h(Y:x,y)$.
For $Z\subset Y$ Zariski-closed   we have 
 $h_v(Y:x,y)=h_v(Y\setminus Z:x,y)+h_v(Z:x,y)$.
For  $f:Z\maps Y$ a Zariski-locally trivial fibre bundle with fibre $F$,  we have 
$h_v(Z:x,y)=h_v(Y:x,y)h_v(F:x,y)$.
Finally $e(Y)=h_v(Y,1,1)$ for any complex variety $Y$.
We denote by     
$$\sum_i (-1)^ib^i_v(Y) y^i=p_v(Y:y):=h_v(Y:y,y)$$  
the  virtual Poincar\'e
polynomial. If $Y$ has pure complex dimension $d$ (or sometimes when $Y$ 
has expected
dimension $d$), we   write
$H_v(Y)=H_v(Y:x,y):=(xy)^{-\frac{d}{2}}h_v(Y:x,y)$, $X_y^v(Y):=H_v(Y:1,y)$ and
$P_v(Y)=P_v(Y:y):=y^{-{d}}p_v(Y:y)$. 
If $Y$ is smooth and projective of dimension $d$ we have therefore
$H_v(Y)=H(Y)$, $X_y^v(Y)=X_y(Y)=\chi_{-y}(Y)$ and $P_v(Y)=P(Y)$.

Let $Y$ be an arbitrary quasiprojective variety (not  necessarily irreducible
or smooth)
over $\C$. We want to show that the Weil conjectures still compute the virtual 
Poincar\'e polynomials. This was pointed out to me by Jun Li, 
and seems to be known
to the experts.

\begin{prop}\label{weilvirt} There  is a finitely generated subring 
 $A=\Z[a_1,\ldots,a_l]\subset \C$ and a variety $Y_A$ over $A$, such that 
$Y=Y_A\times_A \C$, and  the following holds:
 For $m$ a maximal ideal of $A$ we put
$Y_m:=Y_A\times_A A/m$. There is a 
nonempty dense open subset $U$ of
$spec(A)$, such that if $m\in U$ is a maximal ideal of $A$ with quotient field
$\ff_q$, then there exist complex numbers $(a_{i,j})_{i,j}$ with
$|a_{i,j}|=q^{i/2}$, such that for all
$n\in\Z_{>0}$
$$\#Y_q(\ff_{q^n})=\sum_{i}^{}\sum_{j=1}^{b_v^i(Y)} (-1)^i a_{i,j}^n.$$
\end{prop}

\begin{pf}
If $Y$ is smooth and  projective, this is   part of the 
Weil conjectures, proven by Deligne \cite{De}.
The general case is a simple consequence of this and resolution of 
singularities in characteristic $0$. Let $d$ be the largest dimension of a component 
of $Y$.
The proof is by induction on $d$, the case $d=0$ being trivial.
Write $Y=Y_0\sqcup W$, where $Y_0$ is the smooth locus of $Y$, and 
let $\tilde Y=Y_0\sqcup Z$ be a smooth compactification of $Y$.
Then $p_v(Y,z)=p(\tilde Y,z)+p_v(W,z)-p_v(Z,z)$. Let $A=\Z[a_1,\ldots,a_l]\subset \C$ 
be a finitely generated subring, such that $Y$, $\tilde Y$, $Z$, $W$ are already defined
over $A$. Let
$U$ be an open dense subset  of $spec(A)$ where the proposition applies 
to $\tilde Y$ (by the usual Weil conjectures)  $Z$ and  $W$ (by induction). 
Let $m\in U$ be a maximal ideal  with quotient field $\ff_q$.
 Then 
$\#Y_m(\ff_{q^n})=\#\tilde Y_m(\ff_{q^n})+\#W_m(\ff_{q^n})-\#Z_m(\ff_{q^n})$,
and the result follows.
\end{pf}

\subsection{Moduli spaces}
Let again $S$ be an  algebraic surface, 
 $H$ a general ample divisor on $S$, and let $ C\in H^2(X,\Z)$.  
Let $M_S^H(r,C,d)$ denote the moduli space of $H$-semistable sheaves $\E$ on 
$S$ (in the sense of Gieseker-Maruyama), with $c_1(\E)=C$ and discriminant
$d=c_2(\E)-\frac{r-1}{2r}C^2$. Let $M_S^H(r,C,d)_s$ denote the open subspace of
$H$-slope stable sheaves and $N_S^H(r,C,d)$ the subspace of $H$-slope stable
 locally free sheaves.
If $d$ is sufficiently large, then $M_S^H(r,C,d)$ is irreducible and generically
smooth of dimension $e=2rd-(r^2-1)\chi(\oo_S)$ (see e.g. \cite{H-L}).
We put $M_S^H(C,d):=M_S^H(2,C,d)$, $M_S^H(C,d)_s:=M_S^H(2,C,d)_s$ and 
$N_S^H(C,d):=N_S^H(2,C,d)$. 
If $S$ is a rational algebraic surface and $H$ is an ample divisor with 
$HK_S\le 0$, then a slope stable sheaf $\E$ fullfils $\hbox{Ext}^2(\E,\E)=
\hbox{Hom}(\E,\E\otimes K_S)=0$, and therefore  
$M^H_S(r,C,d)_s$ is smooth of dimension 
$e=2rd-(r^2-1)$.

\subsection{Modular forms}

We give a brief review of the results modular forms that we will need.
It might be helpful to also look at \cite{G-Z} section 2.2.
Let $\H:=\big\{\tau\in \C\bigm| \hbox{Im}(\tau)>0\big\}$ be the complex upper
half-plane. For $\tau\in \H$ let $q:=e^{2\pi i \tau}$ and 
$q^{1/n}:=e^{2\pi i \tau/n}$. 
For $a\in \Q$ we often write $(-1)^a$ instead of $e^{\pi i a}$.
We always use the principal
branch of the square root (with $\sqrt{\tau}\in\H$ for $\tau\in \H$ and
$\sqrt{a}\in\R_{>0}$ for
$a\in\R_{>0}$).
We recall the definition of quasimodular forms from \cite{K-Z}.
A modular form of weight
$k$ on a subgroup
$\Gamma\subset Sl(2,\Z)$ of finite index
is a holomorphic function $f$ on $\H$
satisfying
$$f\Big(\frac{a\tau+b}{c\tau+d}\Big)=(c\tau+d)^kf(\tau),\quad \tau\in H, \
\left(\begin{matrix}a&b\\c&d\end{matrix}\right)\in\Gamma$$
growing at most polynomially in $1/\Im(\tau)$ as $\Im(\tau)\to 0$.
An almost holomorphic modular form of weight $k$ is a function 
$F$ on $\H$ with the same transformation properties and growth conditions
as a modular form which is of the form
$F(\tau)=\sum_{m=0}^Mf_m(\tau)(\Im(\tau))^{-m}$ for $M\ge 0$ and 
$f_i$ holomorphic functions. Functions $f$ which occur as  (the holomorphic 
part of $F$)
$f_0(\tau)$ in such an expansion will are called quasimodular forms of weight $k$.
We denote
$\sigma_{k}(n):=\sum_{d|n}d^k$ and by 
$\sigma_{1}^{\text{odd}}(n)$ the sum of the odd divisors of $n$.
For even $k\ge 2$ let 
$$ G_{k}(\tau):=-\frac{B_{k}}{2k}+\sum_{n>0}\sigma_{k-1}(n)q^n
$$ be the  Eisenstein series, where $B_k$ is the $k$-th Bernoulli number. Note
that $G_{k}$ is a modular form of weight $k$ on $SL(2,\Z)$ for $k\ge 4$, but is
only quasimodular for $k=2$, i.e.
$G_2(\tau)+1/(8\pi\Im(\tau))$ is an almost
holomorphic modular form of weight~2. Equivalently
\begin{equation}\label{g2trans}
G_2\bigl(\frac{a\tau+b}{c\tau+d}\bigr)=(c\tau+d)^2G_2(\tau)-\frac{c(c\tau+d)}{4\pi i}
\end{equation}
(see \cite{Z2} p. 242).
Let $\eta(\tau):=q^{1/24}\prod_{n>0}(1-q^n)$ be the Dedekind
eta function and $\Delta:=\eta^{24}$ the discriminant.
We have the transformation laws
\begin{equation}\label{etatrans}
\eta(\tau+1)=(-1)^{1/12}\eta(\tau),\quad
\eta(-1/\tau)=\sqrt{\frac{\tau}{i}}\eta(\tau)\quad \hbox{see \cite{C} VIII.3.}
\end{equation}
We write $y:=e^{2\pi i z}$ for $z$ a complex variable.
Recall the classical theta functions
\begin{align}\label{theta}
\theta_{\mu,\nu}(\tau,z)&:=\sum_{n\in\Z}(-1)^{n\nu}q^{(n+\mu/2)^2/2}
y^{n+\mu/2}
\quad (\mu, \nu\in\{0,1\})
\end{align}
(see e.g. \cite{C} Ch. V, where however the notations and conventions are 
slightly different),
and the   ``Nullwerte"
\begin{align}\label{nullwerte}
\theta(\tau):=\theta_{0,0}(\tau,0)&=\frac{\eta(\tau)^5}
{\eta(\tau/2)^2\eta(2\tau)^2},
\quad
\theta^0_{0,1}(\tau)=\theta_{0,1}(\tau,0)=\frac{\eta(\tau/2)^2}{\eta(\tau)},\\
\theta^0_{1,0}(\tau)=\theta_{1,0}(\tau,0)&=2\frac{\eta(2\tau)^2}{\eta(\tau)},
\quad
\theta_{1,1}(\tau,0)=0.\nonumber
\end{align}
We use the same notations also for $\mu,\nu$ arbitrary in $\Q$.
The  identities (\ref{nullwerte})
follow readily from the product formulas
\begin{align}
\theta_{1,1}(\tau,z)&=
q^{\frac{1}{8}}
(y^{\frac{1}{2}}-y^{-\frac{1}{2}})
\prod_{n>0}(1-q^n)(1-q^ny)(1-q^ny^{-1}),\label{theta11}\\
\theta_{0,1}(\tau,z)
&=\prod_{n>0}(1-q^n)(1-q^{n-\frac{1}{2}}y)(1-q^{n-\frac{1}{2}}y^{-1}),
\nonumber\end{align}
and the fact that $\theta_{\mu,\nu}(\tau,z)= \theta_{\mu,0}(\tau,z+\nu)$.
$\theta_{1,1}$ has the transformation behaviour
\begin{equation}
\label{theta11trans}
\theta_{1,1}(\tau+1,z)=(-1)^{1/4}\theta_{1,1}(\tau,z),\quad 
\theta_{1,1}(-1/\tau, z/\tau)=-i\sqrt{\frac{\tau}{i}}e^{\pi i
z^2/\tau}\theta_{1,1}(\tau,z).
\end{equation}
By the product formulas (\ref{theta11}) we see that
\begin{align}\theta_{0,1}(\tau,z)\theta_{1,1}(\tau,z)&=
\frac{\eta(\tau)^2}{\eta(\tau/2)}\theta_{1,1}(\tau/2,z).
\label{thetahalf}\end{align}
We write 
$$
\widetilde\theta_{1,1}(\tau,z)
:=\frac{\theta_{1,1}(\tau,z)}{y^{\frac{1}{2}}-y^{-\frac{1}{2}}}.
$$
{}From the definitions it is straightforward to see that
\begin{align}\theta_{\mu+2,0}(\tau,z)&=\theta_{\mu,0}(\tau,z),\quad
\theta_{\mu+2,1}(\tau,z)=-\theta_{\mu,1}(\tau,z),\qquad \mu\in
\Q.\label{plus2}\end{align}
By $\frac{(-n+1/2)^2}{4}=\frac{(n-1/2)^2}{4}=(\pm(n/2+1/2))^2$,
one also checks immediately that
\begin{align}\label{theta012}
\theta_{1/2,0}^0(2\tau)&=\sum_{n\in \Z}q^{(n+1/4)^2}=\frac{1}{2}\sum_{n\in \Z}q^{(n+1/2)^2/4}=\frac{\theta_{1,0}^0(\tau/2)}{2}=
\frac{\eta(\tau)^2}{\eta(\tau/2)}.
\end{align}
Following \cite{Go3},\cite{G-Z}, we set 
$  f(\tau):=(-1)^{-1/4}\frac{\eta(\tau)^3}{\theta(\tau)}.
$ 
Let 
$e_2$ and $e_3$ be  
 the $2$-division values of the Weierstra\ss\ $\wp$-function at 
 $\tau/2$ and $(1+\tau)/2$ respectively, i.e.
\begin{align*} 
e_2(\tau)&=\frac{1}{12}+2\sum_{n>0}\sigma^{\text{odd}}_1(n)q^{n/2},\\
e_3(\tau)&=\frac{1}{12}+2\sum_{n>0}(-1)^n\sigma^{\text{odd}}_1(n)q^{n/2},
\end{align*} 
(see e.g. \cite{H-B-J} p 132). It is easy to see that 
$e_3(2\tau+1)=e_2(2\tau)$. We also see that $\theta(2\tau+1)=\theta_{0,1}^0(2\tau)$ 
and 
$f(2\tau+1)=\eta(2\tau)^4/\eta(\tau)^2$.
We write  
\begin{align}
\label{utau} u(\tau)&:=-\frac{f(\tau)^2}{3e_3(\tau)}, \quad
\overline
u(\tau):=u(2\tau+1)=-\frac{\eta(2\tau)^8}{3e_2(2\tau)\eta(\tau)^4}.\end{align}

\begin{rem} 
\label{gammau} Let 
\begin{equation*}
T:=\left(\begin{matrix}1&1\\0&1\end{matrix}\right),\quad
V:=T^2=\left(\begin{matrix}1&2\\0&1\end{matrix}\right),\quad
S:=\left(\begin{matrix}0&-1\\1&0\end{matrix}\right).
\end{equation*}
Let $\Gamma_u=\pm\<V^2,VS,SV\>$; this is a subgroup of index $6$ of
$SL(2,\Z)$. $u(\tau)$ is a modular function on
$\Gamma_u$. Let 
$\Gamma(2):=\big\{A\in Sl(2,\Z)\bigm| A\equiv id \hbox{ mod } 2\big\}.$
Let $
X:=\Big(\begin{matrix}2&1\\0&1\end{matrix}\Big).
$
It is easy to see that $X^{-1}\Gamma_u X=\Gamma(2)$. In other words 
a function $g(\tau)$ is a modular function on $\Gamma_u$, if and only if 
$h(\tau):=g(2\tau+1)$ is a modular function  on $\Gamma(2)$.
In particular $\overline u(\tau)$ is a modular function on $\Gamma(2)$.
\end{rem}

\subsection{Theta functions for indefinite lattices}

We review the definition of theta functions for indefinite lattices 
from \cite{G-Z}.
Let $\Gamma$ be a lattice, i.e. a free $\Z$ module $\Gamma $ together with
an $\Z$-valued bilinear form  $\<x,y\>$ on $\Gamma$.
The extension of the bilinear form  to 
$\Gamma_\C:=\Gamma\otimes\C$ and $\Gamma_{\R}=\Gamma\otimes \R$ 
is  denoted in the same way.
The {\it type} of $\Gamma$ is the pair $(r-s,s)$, where $r$ is the
rank of $\Gamma$ and $s$ the 
largest rank of a sublattice of $\Gamma$ on which $\<\ ,\ \>$ is negative definite.
Let $M_\Gamma$ be the space of meromorphic functions on $\H\times \Gamma_\C$.
For $v\in \Gamma_\Q$, $A=\Big(\begin{matrix}a&b\\c&d\end{matrix}\Big)$,
and $k\in\Z$ we
put
\begin{align}\label{strich}
f|v(\tau,x)&:=q^{\<v,v\>/2}\exp(2\pi i \<v,x\>)f(\tau,x+v\tau),\\
f|_kA(\tau,x)&:=(c\tau+d)^{-k}\exp\Big(-\pi i \frac{\<x,x\>}{c\tau+d}\Big)
f\Big(\frac{a\tau+b}{c\tau+d},\frac{x}{c\tau+d}\Big).
\end{align}
Now assume that $\Gamma$ is unimodular of type $(r-1,1)$.   We fix a vector 
$f_0\in \Gamma_\R$ with $\<f_0, f_0\><0$, and let  
\begin{align*} C_\Gamma&:=\big\{f\in \Gamma_{\R}\bigm|\<f,f\><0,\;\,\<f,
f_0\><0\big\},\\
 S_\Gamma&:=\big\{f\in \Gamma\bigm| 
f \hbox{ primitive}, \ \<f,f\>=0,\ \<f,f_0\><0\big\}. 
\end{align*}
For $f\in S_\Gamma$ put
$$
D(f):=
\big\{(\tau,x)\in \H\times \Gamma_\C\bigm| 0<\Im(\<f, x\>)<\Im(\tau) \ \big\},
$$
and for $f\in C_\Gamma$ put $D(f):=\H\times \Gamma_\C$.
For $t\in \R$ we put 
$\mu(t):=1$,  if $t\ge 0$ and $\mu(t)=0$ otherwise.
Let $f,g\in C_\Gamma\cup S_\Gamma$. For 
$c\in \Gamma$ and $(\tau,x)\in D(f)\cap D(g)$ 
we put 
\begin{equation}
\label{thetafgcb}
\Theta^{f,g}_{\Gamma,c}(\tau,x):=
\sum_{\xi\in \Gamma+c/2}\bigl(\mu(\<\xi, f\>)
-\mu(\<\xi, g\>)\bigr)\,q^{\<\xi,\xi\>/2}\,e^{2\pi i \<\xi, x\>},
\end{equation}
and 
$\Theta^{f,g}_{\Gamma}:=\Theta^{f,g}_{\Gamma,0}$.

Assume now that  $f,g\in S_\Gamma$. Then (see \cite{G-Z})
 the function $\Theta^{f,g}_{\Gamma,c,b}$ has a meromorphic extension to $\H\times
\Gamma_\C$, which is defined as follows. 
Let 
$$F:\H\times\C^2\to \C;\ (\tau,u,v)\mapsto 
\frac{\eta(\tau)^3\theta_{1,1}(\tau,(u+v)/(2\pi i))}
{\theta_{1,1}(\tau,u/(2\pi i))\theta_{1,1}(\tau,v/(2\pi i))},$$
(see \cite{Z1}; note the different conventions for $\theta_{1,1}$ in \cite{Z1}).
We have
$$F(\tau,{u},{v})=
\sum_{n\ge 0, m>0} q^{nm}e^{-nu-mv}-
\sum_{n> 0, m\ge 0} q^{nm}e^{nu+mv},$$
(see \cite{G-Z} section 3.1).
Assume $\<f,g\>=-N\in\Z_{<0}.$
We denote by $[f,g]$ the lattice generated by $f$ and $g$ and by $[f,g]^{\perp}$
its orthogonal complement. Let $L:=[f,g]\oplus [f,g]^{\perp}$.
For $(\tau,x)\in \H\times \Gamma_\C$,
we put
\begin{equation}
\Theta^{f,g}_L(\tau,x):=F(N\tau,-2\pi i\<f,x\>,2\pi i\<g,x\>) \Big(\sum_{\xi \in
[f,g]^{\perp}}q^{\<\xi,\xi\>/2}e^{2\pi i \<\xi,x\>}\Big).\label{thetaext1}
\end{equation}
Let $P$ be a system of representatives of $\Gamma$ modulo $L$.
Then, using the notation of (\ref{strich}), 
the meromorphic extension is given by 
\begin{equation}\label{thetaext}
\Theta^{f,g}_{\Gamma,c}:=\sum_{t\in P}\Theta^{f,g}_L(\tau,x)|(t+c/2).
\end{equation}
For $f,g\in S_\Gamma$ the following is shown in  \cite{G-Z}:
For $|\Im(\<f, x\>)/\Im(\tau)|<1$, $|\Im(\<g, x\>)/\Im(\tau)|<1$ we have
\begin{align*}
\Theta^{f,g}_{\Gamma,c}(\tau,x)&=
\frac{1}{1- e^{2\pi i\<f,x\>}}\sum_{{ \<\xi, f\>=0}\atop
{\<f, g\>\le  \<\xi, g\><0}} q^{ \<\xi,\xi\>/2}e^{2\pi i \<\xi,x\>}
-\frac{1}{1- e^{ 2\pi i \<g,x\>}}\sum_{{ \<\xi,  g\>=0}\atop
{\<f, g\>\le  \<\xi, f\><0}} q^{ \<\xi,\xi\>/2}e^{ 2\pi i\<\xi,x\>}\\
&\qquad +\sum_{ \xi\cdot  f> 0> \xi\cdot g}
 q^{ \<\xi,\xi\>/2}\big(e^{2\pi i\<\xi, x\>}-e^{-2\pi i \<\xi,x\>}\big).
\end{align*}
Here the sum is taken over all $\xi\in \Gamma+c/2$.
For $b,c\in \Gamma$ and any characteristic vector $w$ of $\Gamma$
we have
\begin{align}{ }
(\Theta_{\Gamma,c}^{f,g}/\theta^{\sigma(\Gamma)})|S(\tau,x+b/2)
&=(\Theta_{\Gamma,b}^{f,g}/\theta^{\sigma(\Gamma)})(\tau,x+c/2),\nonumber\\
\label{thetatrans}\Theta_{\Gamma,c}^{f,g}(\tau+1,x)
&=(-1)^{3\<c,c\>/4-\<c,w\>/2}
\Theta_{\Gamma,c}^{f,g}(\tau,x+(w-c)/2),\\
\Theta_{\Gamma,c}^{f,g}(\tau+2,x)&=(-1)^{\<c,c\>/2}
\Theta_{\Gamma,c}^{f,g}(\tau,x).
\nonumber
\end{align}
The last two formulas are elementary consequences of the 
definition (\ref{thetafgcb}),
which also hold for $f,g\in C_\Gamma\cup S_\Gamma$.

\subsection{Hilbert schemes}
For a general algebraic surface $S$, we
 denote by $S^{[n]}$ the Hilbert scheme of subschemes of 
length $n$ on $S$. $S^{[n]}$ is smooth of dimension $2n$ \cite{F}, and
its Hodge numbers have been computed (\cite{E-S}, \cite{Go1}, \cite{G-S},
\cite{Ch}). Using (\ref{theta11}),
the results can be easily translated to
\begin{align}\label{hilbhodge}
\sum_{n\ge 0}X_y(S^{[n]})q^{n-e(S)/24}&=
\frac{\eta(\tau)^{\sigma(S)-\chi(\oo_S)}}
{\tilde\theta_{1,1}(\tau,z)^{\chi(\oo_S)}}.
\end{align}
(Recall that we write $y:=e^{2\pi i z}$).
In particular 
$$\sum_{n\ge 0}e(S^{[n]})q^{n-e(S)/24}=\frac{1}{\eta(\tau)^{e(S)}},\qquad
\sum_{n\ge 0}\sigma(S^{[n]})(-1)^nq^{n-e(S)/24}=\frac{\eta(\tau)^{\sigma(S)}}
{\eta(2\tau)^{2\chi(\oo_S)}}.$$

\section{Relation to locally free sheaves and blowup formulas}

In this section let $S$ be an arbitrary smooth projective surface, and let $C\in
H^2(S,\Z)$. Let $\widehat S$ be the blowup of $S$ in a point and $E$ the 
exceptional divisor. Let $H$ be a general ample divisor on $S$ (general means that
it does not lie on a wall with respect to $(r,C)$, see \cite{Y3}; in the 
case $r=2$ we will discuss walls and chambers in the next section).
 We will usually denote the cohomology classes
on $S$ and their pullbacks to $\widehat S$ by the same letter. 
We denote by $M_{\widehat S}^H(r,C+bE,d)_s$ the space of slope stable
sheaves on $\widehat S$ which are stable with respect to (the pullback of)
$H$. It  can be identified with 
$M_{\widehat S}^{H-\epsilon E}(r,C+bE,d)_s$ for $\epsilon >0$ small enough.
 
We want to relate the virtual Poincar\'e polynomials of 
$M_S^H(r,C,d)_s$, $N_S^H(r,C,d)$ and $M_{\widehat S}^H(r,C+bE,d)_s$.
In fact we will see that the generating function for $\widehat S$ is obtained
from that for $S$ by multiplying by a suitable theta function 
and dividing by a power of the eta function.  
The results are easy consequences of corresponding results of 
Yoshioka about the counting of points of these moduli spaces over finite fields
and of Prop.~\ref{weilvirt}.
We write
$$P_v(M^H_S(r,C,d)_s)=y^{-e} p_v(M^H_S(r,C,d)_s,y), \ 
P_v(N^H_S(r,C,d))=y^{-e} p_v(N^H_S(r,C,d),y), $$
where 
$e=2rd-(r^2-1)\chi(\oo_S)$ is the virtual dimension,
which agrees with the actual dimension for $d$ sufficiently large.

\begin{prop}\label{blowprop}
Let $S$ be an algebraic surface and let $H$ be a general ample divisor on $S$.
\begin{enumerate}
\item
$$\displaystyle{
\sum_{d\ge 0}P_v(M^H_S(r,C,d)_s)q^d
=\left(\prod_{k\ge
1}\prod_{b=1}^r\prod_{i=0}^4 (1-y^{i-2b}q^k)^{(-1)^{i+1}b_i(S)}
\right)\left(\sum_{d\ge 0}P_v(N^H_S(r,C,d))q^d\right),
}$$
 in particular
$$\sum_{d\ge 0}e(M^H_S(r,C,d)_s)q^d=\frac{q^{re(S)/24}}{\eta(\tau)^{re(S)}}
\left(\sum_{d\ge 0}e(N^H_S(r,C,d))q^d\right).$$
\item
Let $A=(a_{ij})_{ij}$ be the $(r-1)\times (r-1)$-matrix
with  entries $a_{ij}=1$ for $i\le  j$ and 
$a_{ij}=0$  otherwise. We view 
elements of $\R^{r-1}$   as column vectors.
We write ${I}$ for the column vector of length $r-1$ with all entries equal to
one. Then 
$$
\sum_{d\ge 0}P_v(M^H_{\widehat S}(r,C+bE,d)_s)q^d
=
\frac{q^{{r}/{24}}}{\eta(\tau)^{r}}
\left(\sum_{v\in \Z^{r-1}+\frac{b}{r}I}(y^2)^{v^t A I} q^{v^t A v}\right)
\left(\sum_{d\ge 0}P_v(M^H_{S}(r,C,d)_s)q^d\right),
$$
in particular 
$$
\sum_{d\ge 0}e(M^H_{\widehat S}(r,C+bE,d)_s)q^d
=\frac{q^{{r}/{24}}}{\eta(\tau)^{r}}\left(\sum_{v\in \Z^{r-1}+\frac{b}{r}I} 
q^{v^t Av}\right)\left(\sum_{d\ge 0}e(M^H_{S}(r,C+bE,d)_s)q^d\right).
$$
\end{enumerate}
\end{prop}

\begin{pf}
(1) is a consequence of (\cite{Y1}, Thm.~0.4) and  Prop.~\ref{weilvirt}:
Let $X$ be a surface over $\ff_q$.
For every sheaf $E$ in $M^H_X(r,C,d)_s(\ff_q)$
there is an exact sequence
$0\to E\to E^{\vee\vee}\to E^{\vee\vee}/E\to 0$, where $E^{\vee\vee}\in 
N^H_S(r,C,d-k)_s(\ff_q)$ and $E^{\vee\vee}/E\in \Quot^k_{E^{\vee\vee}}(\ff_q)$
for a suitable $k\le d$.
In fact it is easy to see that if $E$ is defined over $\overline \ff_q$,
then it is defined over $\ff_q$ if and only if both $E^{\vee\vee}$ and 
$E^{\vee\vee}/E$ are. 
For a   sheaf $F$ over $X$ we denote by $\Quot^k_{F}$ the (Grothendieck)
scheme of 
quotients of length $k$ of $F$ and by  $\Quot^k_{F,p}$ the subscheme
 (with the reduced structure) of quotients with support in the point $p\in X$.
If $F$ is locally free of rank $r$ and $p$ is defined over $\ff_{q}$, we get
isomorphisms $\Quot^k_{F,p}\simeq \Quot^k_{\oo^{\oplus r}_X,p}$ over $\ff_{q}$. 
In particular $\#\Quot^k_{F,p}(\ff_{q})=\#\Quot^k_{\oo^{\oplus r}_X,p}(\ff_{q}).$
Therefore the proof of (\cite{Y1}, Thm.~0.4) for the numbers
$\#\Quot^k_{\oo^{\oplus r}_X}(\ff_{q})$ can be repeated for 
$\#\Quot^k_{F}(\ff_{q})$, the only numbers entering the calculation 
being the $\#\Quot^k_{F,p}(\ff_{q^n})$. Therefore 
$\#\Quot^k_{F}(\ff_{q})=\#\Quot^k_{\oo^{\oplus r}_X}(\ff_{q})$
(see also Y1, p.194).
This gives 
$$\# M^H_X(r,C,d)_s(\ff_q)=\sum_{k\le d} \# N^H_X(r,C,d-k)_s(\ff_q)\cdot
\# \Quot^k_{\oo^{\oplus r}_X}(\ff_q).$$
Applying Prop.~\ref{weilvirt} to a good reduction $X$ of $S$ modulo $q$,
we obtain immediately
$$\sum_{d\ge 0} \sum_{d\ge 0}p_v(M^H_S(r,C,d)_s)q^d
=\left(\prod_{k\ge
1}\prod_{b=1}^r\prod_{i=0}^4 (1-y^{2rk+i-2b}q^k)^{(-1)^{i+1}b_i(S)}
\right)\left(\sum_{d\ge 0}p_v(N^H_S(r,C,d))q^d\right),$$
(recall the signs in the definition of $p_v$). By the definition
of $P_v$ and the formula $e=2rd-(r^2-1)\chi(\oo_S)$,
we see  that in order to replace 
$p_v$ by $P_v$ we have to replace the factor $(1-y^{2rk+i-2b}q^k)$
by $(1-y^{i-2b}q^k)$.

(2) We apply Prop.~\ref{weilvirt} to  (\cite{Y3}, Prop. 3.4).
Using again $e=2rd-(r^2-1)\chi(\oo_S)$  we obtain 
\begin{align*}
\sum_{d\ge 0}&P_v(M^H_{\widehat
S}(r,C+bE,d)_s)q^d\\&=\frac{q^{{r}/{24}}}{\eta(\tau)^{r}}
\left(\sum_{(a_1,\ldots,a_{r})} (y^2)^{w(a_1,\ldots,a_{r})} 
q^{-\sum_{i<j}a_ia_j}
\right)\left(\sum_{d\ge 0}P_v(M^H_{S}(r,C,d)_s)q^d\right).
\end{align*}
Here the sum runs through the $r$-tuples $(a_1,\ldots,a_{r})$ in
$\Z+\frac{b}{r}$ with $\sum_{i=1}^r a_i=0$,
and 
$$w(a_1,\ldots,a_r)=\sum_{i<j\le r}\binom{a_j-a_i}{2}+r{\sum_{i<j\le
r}a_ia_j}.$$ We note that equivalently we can let the sum run through the 
$(r-1)$-tuples $(a_1,\ldots,a_{r-1})$, and put $a_r=-\sum_{i=1}^{r-1}a_i$.
Then 
$$
{-\sum_{i<j\le r}a_ia_j}=\sum_{j\le i\le r-1}a_ia_j.
$$
Furthermore we have 
$$
\sum_{i<j\le r}(a_j-a_i)^2=2r\sum_{j\le i \le r-1}a_ia_j
$$
and 
$$
\sum_{i<j\le r}(a_j-a_i)=-2\left(\sum_{i=1}^{r-1}(r-i)a_i\right).
$$
Putting things together, we obtain
$$
w(a_1,\ldots, a_{r})=\sum_{i=1}^{r-1}(r-i) a_i=(a_1,\ldots,a_{r-1})A I.
$$
Finally we note that 
$$
\sum_{j\le i\le r-1} a_ia_j=(a_1,\ldots,a_{r-1})A(a_1,\ldots,a_{r-1})^t.
$$
\end{pf}

\begin{rem}\begin{enumerate}
\item
Li and Qin (\cite{L-Q1}, \cite{L-Q2}) have shown a blowup formula 
for the virtual Hodge
polynomials in the case $r=2$ using completely different methods. 
In particular they
also obtain
a blowup formula for the Euler numbers. Their method also gives
 a blowup formula for the
virtual Hodge polynomials of the Uhlenbeck compactification.
We write again
$H_v(M^H_S(r,C,d)_s)=(xy)^{-e/2} h_v(M^H_S(r,C,d)_s)$ with 
$e=2rd-(r^2-1)\chi(\oo_S)$.
Then, writing $x=e^{2\pi i u}$,  their result can be rewritten as 
\begin{align*}
\sum_{d\ge 0}H_v(M^H_{\widehat
S}(C,d)_s)q^d&=\frac{q^{1/12}\theta_{0,0}(2\tau,u+z)}{\eta(\tau)^2}
\left(\sum_{d\ge 0}H_v(M^H_S(C,d)_s)q^d\right),\\
\sum_{d\ge 0}H_v(M^H_{\widehat
S}(C+E,d)_s)q^d&=\frac{q^{1/12}\theta_{1,0}(2\tau,u+z)}{\eta(\tau)^2}
\left(\sum_{d\ge 0}H_v(M^H_S(C,d)_s)q^d\right).
\end{align*}
This is the case $r=2$ of the formula 
\begin{align*}
\sum_{d\ge 0}&H_v(M^H_{\widehat S}(r,C+bE,d)_s)q^d\\&=
\frac{q^{{r}/{24}}}{\eta(\tau)^{r}}
\left(\sum_{v\in \Z^{r-1}+\frac{b}{r}I}(xy)^{v^t AI} q^{v^t A v}\right)
\left(\sum_{d\ge 0}H_v(M^H_{S}(r,C,d)_s)q^d\right).
\end{align*}
I expect that this formula holds for all $r$.
\item
Using \cite{Y5},  Prop.~\ref{blowprop}(2)
can also be rewritten:
Let $A_{r-1}=\big\{(x_1,\dots, x_r)\bigm| \sum_i x_i=0\big\}$ be the
$A_{r-1}$-lattice and
$e_1,
\dots, e_{r-1}$ its standard basis.  Let
$a:=\sum_{i=1}^{r-1} i(r-i)
e_i$ and
$\lambda=(1-1/r,-1/r,\dots,-1/r)$.
Then the theta function on the left hand side in Prop. \ref{blowprop} can be
written as
$$
\sum_{v \in A_{r-1}+b \lambda}y^{\<v,a\>}q^{\<v,v\>/2},
$$
where $\<\ ,\  \>$ is the pairing of $A_{r-1}$.
This was pointed out to me by K.~Yoshioka.
\end{enumerate}
\end{rem}

\section{Wallcrossing and theta functions}

\subsection{Wallcrossing}

Now let $S$ again be a rational algebraic surface. Let   $\Gamma$ be the lattice 
$H^2(S,\Z)$
with the negative of the intersection form as quadratic form, i.e. 
for $A,B\in \Gamma$ let $\<A,B\>=-AB$.
In this section we want to relate the Hodge numbers of the 
moduli spaces $M^H_S(C,d)$ to the theta functions
$\Theta_{\Gamma,C}^{F,H}$ from \cite{G-Z}.
The dependence of the moduli spaces $M^H_S(C,d)$ on the polarization $H$
and the corresponding dependence of the Donaldson invariants
 has been studied by a number of authors
\cite{Q1}, \cite{Q2}, \cite{F-Q}, \cite{Go2}, \cite{E-G}, \cite{Y3}, \cite{L}.
We follow (with some modifications) the notations in \cite{Go2}, \cite{E-G}.

An ample divisor $H$ is called {\it good} if $K_S\cdot H\le 0$.
We denote by $\cc_S$ the ample cone of $S$ and by 
$\cc_S^G$ the subcone of all good ample divisors.
A class $\xi\in H^2(X,\Z)+C/2$ is called  of type
$(C,d)$ if $\xi^2+d\in \Z_{\ge 0}$.
In this case we call $W^\xi:=\xi^\perp\cap \cc_S$ the  
{\it wall} defined by $\xi$. 
If $\xi^\perp\cap \cc_S^G\ne \emptyset$, we call $W^\xi$ a good wall.
The chambers of type
$(C,d)$ are the connected  components of the complement of the walls of 
type $(C,d)$ in
$\cc_S$.  If $L$ and $H$ lie in the same chamber of type $(C,d)$, then
$M_S^L(C,d)=M_S^H(C,d)$.
We say that $L$ lies on a wall of type $C$, if $L\xi=0$ for some class
$\xi\in H^2(X,\Z)+C/2$.

\begin{thm}\label{walltheta} Let $C\in H^2(S,\Z)$.
Let $H,L\in \cc^G_S$   not   on a wall of type $C$. 
Then
\begin{enumerate}
\item
\begin{align*}
\sum_{d\ge 0}\big( X^v_y(M_S^H(C&,d))-X^v_y(M_S^L(C,d))\big)q^{d-e(S)/12}
\\&=\frac{
\eta(\tau)^{2\sigma(S)-2}(y^{\frac{1}{2}}-y^{-\frac{1}{2}})}
{\theta_{1,1}(\tau,z)^2}\Theta^{L,H}_{\Gamma,C}(2\tau,K_Sz),\\
\sum_{d\ge 0}
(e(M^H_S(C&,d))-e(M^L(C,d)))q^{d-e(S)/12} =\frac{1}{\eta(\tau)^{2e(S)}}
\text{\rm Coeff}_{2\pi i z}\big(\Theta^{L,H}_{\Gamma,C}(2\tau, K_Sz)\big).
\end{align*}

\item Assume now that $C\not \in 2H^2(S,\Z)$. Then  we can replace 
$X_y^v$ by $X_y$ in (1). Furthermore 
\begin{align*}
\sum_{d\ge 0} (e(N^H_S(C,d))-e(N^L_S(C,d)))q^d &=\text{\rm Coeff}_{2\pi i
z}\big(\Theta^{L,H}_{\Gamma,C}(2\tau, K_Sz)\big)\\
\sum_{d\ge
0}(-1)^{e(d)/2}(\sigma(M_S^{H}(C,d))-\sigma(M_S^{L}(C,d)))q^{d-e(S)/12}&=
\frac{\eta(\tau)^{2\sigma(S)}}{2i\eta(2\tau)^{4}}
\Theta^{L,H}_{\Gamma,C}(2\tau, K_S/2).
\end{align*} Here $e(d):=4d-3$ is the dimension of
$M_S^{H}(C,d)$. 
\end{enumerate}
\end{thm}

\begin{pf}

This is essentially a reformulation of Thm.~3.4 from \cite{Go2}.
 Assume that 
$H$ and $L$ do not lie on a wall of type $C$. 
The  result of \cite{Go2} gives 
$$
y^{2d-3/2}(X_y^v(M_S^{H}(C,d))-X^v_{y}(M_S^{L}(C,d)))
=\sum_{\xi} y^{d+\xi^2}X_{y}((S\sqcup S)^{[d+\xi^2]}) y^{d-\xi^2}
\frac{y^{\xi K_S}-y^{-\xi K_S}}{y(y-1)},
$$
where the sum runs through all classes of type $(C,d)$ with
$\xi H<0<\xi L$. We sum over all $d\ge 0$.
We use (\ref{hilbhodge}), noting that  
$\sum_{n\ge 0}X_y((S\sqcup S)^{[n]})q^n
=\big(\sum_{n\ge 0}X_y(S^{[n]})q^n\big)^2.$
We obtain
\begin{equation}
\sum_{d\ge 0}\big(X^v_y(M_S^{H}(C,d))-X^v_y(M_S^{L}(C,d))\big)
q^{d-e(S)/12}=\frac{\eta(\tau)^{2\sigma(S)-2}}
{\tilde\theta_{1,1}(\tau,z)^2}
\left(\sum_{\xi }q^{-\xi^2}
\frac{y^{\xi K_S}-y^{-\xi K_S}}
{y^{\frac{1}{2}}-y^{-\frac{1}{2}}}\right).\label{wallform}
\end{equation}
The sum on the right hand side  runs
 through all $\xi\in H^2(X,\Z)+C/2$ satisfying
$\xi H<0<\xi L$.
Using the definition (\ref{thetafgcb}) of the theta functions 
$\Theta^{g,f}_{\Gamma,c}$,
we obtain.  
\begin{align*}
\frac{\theta_{1,1}(\tau,z)^2}
{\eta(\tau)^{2\sigma(S)-2}
(y^{\frac{1}{2}}-y^{-\frac{1}{2}})}&
\left(\sum_{d\ge 0}\big( X^v_y(M_S^H(C,d))-X^v_y(M_S^L(C,d))\big)
q^{d-e(S)/12}\right)\\
&=\sum\Sb\xi\in H^2(S,\Z)+C/2\\
\xi H<0<\xi L\endSb q^{-\xi^2}\big(y^{\xi K_S}-y^{-\xi K_S}\big)
\\
&=\sum\Sb \xi\in \Gamma+C/2\\ \<\xi,H\><0<\<\xi,L\>\endSb
q^{\<\xi,\xi\>}\big(y^{\<\xi,K_S\>}-y^{-\<\xi,K_S\>}\big)
\\
&=\sum_{\xi\in \Gamma+C/2}
q^{\<\xi,\xi\>}(\mu(\<\xi,L\>)-\mu(\<\xi,H\>))y^{\<\xi,K_S\>}
\\
&=\Theta^{L,H}_{\Gamma,C}(2\tau,K_S z).
\end{align*}
 
Specializing (\ref{wallform}) to the Euler number we obtain:
\begin{align*}
\sum_{d\ge 0}\big(e(M_S^{H}&(C,d))-e(M_S^{L}(C,d))\big)
q^{d-e(S)/12}=
\frac{1}{\eta(\tau)^{2e(S)}}
\left(\sum_{H\xi<0<L\xi }2\xi K_S q^{-\xi^2}\right)\\
&=\frac{1}{\eta(\tau)^{2e(S)}}
\left(\sum_{\<H,\xi\><0<\<L,\xi\>}2\<\xi,K_S\> q^{\<\xi,\xi\>}\right)\\
&=\frac{1}{\eta(\tau)^{2e(S)}}
\hbox{\rm Coeff}_{2\pi i z}(\Theta^{L,H}_{\Gamma,C}(2\tau,K_S z))
\end{align*}
The last line follows directly from  (\ref{thetafgcb}).

Assume now that $C\not \in 2H^2(S,\Z)$. Then $M^H_S(C,d)=M^H_S(C,d)_s$ and  
$M^L_S(C,d)=M^L_S(C,d)_s$
are smooth, and we can replace $X^v_y(M^H_S(C,d))$ by $X_y(M^H_S(C,d))$
and $X^v_y(M^L_S(C,d))$ by $X_y(M^L_S(C,d))$. 
Furthermore we get by Prop.~\ref{blowprop}
$$\sum_{d\ge 0} (e(N^H_S(C,d))-e(N^L_S(C,d)))q^d
=\hbox{\rm Coeff}_{2\pi i z}\big(\Theta^{L,H}_{\Gamma,C}(2\tau, K_Sz)\big).
$$
Finally we obtain from  (\ref{wallform})
\begin{align*}
\sum_{d\ge 0}(-1)^{-e(d)/2}(\sigma(M_S^{H}&(C,d))-\sigma(M_S^{L}(C,d)))q^{d-e(S)/12}\\
&=
\frac{\eta(\tau)^{2\sigma(S)}}{ \eta(2\tau)^{4}}
\left(\sum_{\xi H<0<\xi L} q^{-\xi^2}\frac{(-1)^{\xi K_S}-(-1)^{-\xi
K_S}}{(-1)^{1/2}-(-1)^{-1/2}}\right)\\
&=\frac{\eta(\tau)^{2\sigma(S)}}{2i\eta(2\tau)^{4}}
\left(\sum_{\<\xi,H\><0<\<\xi,L\>}
q^{\<\xi,\xi\>}\big((-1)^{\<\xi,K_S\>}-(-1)^{-\<\xi,K_S\>}\big)\right)\\
&=\frac{\eta(\tau)^{2\sigma(S)}}{2i\eta(2\tau)^{4}}
\Theta^{L,H}_{\Gamma,C}(2\tau,K_S/2).
\end{align*}
As the signature can only be nonzero if $e(d)$ is even, we can replace $(-1)^{-e(d)/2}$
by $(-1)^{e(d)/2}$.
\end{pf}

\subsection{Extension of the invariants}
A class $F\in H^2(X,\Z)$ is called nef if its intersection with every effective curve
is nonnegative. The real cone $\overline \cc_{S,\R}$ of nef classes is the closure
of the ample cone. Let  $\delta(\cc_S)$ be the set of primitive
classes in
$(\overline \cc_{S,\R}\setminus \cc_{S,\R})\cap H^2(S,\Z)$, and put
$\overline \cc_S:=\cc_S\cup\delta(\cc_S)$. Let 
$$\SS_S:=\big\{ F\in \delta( \cc_S)\bigm| F^2= 0\big\}.$$
An upper index $G$ will indicate
that we allow only classes $H$ with $HK_S\le 0$.
We   now  extend the generating functions for  the $\chi_y$-genera,
Euler numbers and signatures of the $M^H_S(C,d)$ to the whole of
$\overline \cc_S$.
 
\begin{defn}\label{extend} Let $C\in H^2(S,\Z)$.
\begin{enumerate}
\item 
Let $F \in \delta(\cc^G_S)$, and assume that $CF$
is odd.  Then, for each $d$, the class $F$ lies in the closure of a unique
chamber $\alpha\subset \cc_S^G$ of type $(C,d)$. We put $M^F_S(C,d):=M^H_S(C,d)$ 
for $H\in \alpha$. 
\item
Let $C\in H^2(S,\Z)$, and let $H\in \cc_S^G$, not lying on a wall of type $C$.
Let  
$F\in \overline \cc_S$. If $F\in \delta(\cc_S)$, we assume that $K_SF\ne 0$ or 
$FC$ is odd. We put
$$
\xx^{S,F}_C:=
\left(\sum_{d\ge
0}X^v_y(M^H_S(C,d))q^{d-e(S)/12}\right)+\Theta^{H,F}_{\Gamma,C}(2\tau,K_S z)
\frac{
\eta(\tau)^{2\sigma(S)-2}(y^{\frac{1}{2}}-y^{-\frac{1}{2}})}
{\theta_{1,1}(\tau,z)^2}.
$$
If $C\not \in 2H^2(S,\Z)$ we can replace $X^v_y(M^H_S(C,d))$ by 
$X_y(M^H_S(C,d))$.
We also put 
\begin{align*}(\ee^{S,F}_C)_0&:=\sum_{d\ge 0}e(N^H_S(C,d))q^d+
\hbox{\rm Coeff}_{2\pi i z} 
\big(\Theta^{H,F}_{\Gamma,C}(2\tau,K_Sz)\big),\qquad
\ee^{S,F}_C:=\frac{(\ee^{S,F}_C)_0}{\eta(\tau)^{2e(S)}}
\end{align*}
Finally we put $\Sigma^{S,F}_C=0$ if $C^2\equiv 0$ modulo $2$, and otherwise
\begin{align*}
\Sigma^{S,F}_C:=
&\sum_{d\ge
0}(-1)^{e(d)/2}\sigma(M_S^{H}(C,d))q^{d-e(S)/12}
+\frac{\eta(\tau)^{2\sigma(S)}}{2i\eta(2\tau)^{4}}
\Theta^{H,F}_{\Gamma,C}(2\tau,K_S/2).
\end{align*}
The cocycle condition
\begin{align}\label{cocycle}
\Theta^{F,G}_{\Gamma,C}+\Theta^{G,H}_{\Gamma,C}&=\Theta^{F,H}_{\Gamma,C}
\qquad\hbox{(Rem.~3.4 from \cite{G-Z})}\end{align}
 and Thm.~\ref{walltheta} imply that the definitions
of $\xx^{S,F}_C$, $(\ee^{S,F}_C)_0$, $\ee^{S,F}_C$ and $\Sigma^{S,F}_C$
are independent of $H$.
\end{enumerate}
\end{defn}

\begin{rem}
The definitions above are motivated as follows.
Denote by $X^F_S(C,d)$, $E_S^F(C,d)$ and $S_S^F(C,d)$ 
the coefficients of  $q^{d-e(S)/12}$ in 
$ \xx^{S,F}_C$, $\ee^{S,F}_C$  and  $\Sigma^{S,F}_C$.
Then 
\begin{enumerate}
\item If $F\in \cc_S^G$ does not lie on a wall of type $(C,d)$, then
by Thm.~\ref{walltheta}
$X^F_S(C,d)=X_y(M^F_S(C,d))$,  
$E_S^F(C,d)=e(M^F_S(C,d))$ and $S_S^F(C,d)=(-1)^{e(d)/2} \sigma(M^F_S(C,d))$.
If $C\equiv 0$ modulo $2$, then the moduli space $M^F_S(C,d)$ will sometimes
be singular. We defined $\sigma(M^F_S(C,d)):=0$, which is reasonable, as
$M^F_S(C,d)$ has odd complex dimension.
\item If $L\in \cc_S$ fulfills $K_S L>0$, then $M^L_S(C,d)$ is not 
necessarily smooth,
and the wallcrossing formulas of Thm.~\ref{walltheta} need not be true. 
We extend the
generating function by formally requiring that  Thm.~\ref{walltheta} holds. 
We have seen
above that this gives a consistent definition of the generating functions. 
I believe that there exists a   geometric definition    of the $X_y(M_S^L(C,d))$ 
given by this generating
function.

\item If $F\in \cc_S$ lies on a finite number of walls of type $(C,d)$,
then  $X^F_S(C,d)$ is the average of the $X_y(M^G_S(C,d))$ over all the
chambers of type $(C,d)$ which contain $F$ in their closure
(and similarly for  $E^F_S(C,d)$ and $S^F_S(C,d)$).
\item If $F\in \SS_S$ and $CF$ is even, then $F$ usually lies on infinitely
many walls of type $(C,d)$ and in the closure of infinitely many chambers.
We can view $X_S^F(C,d)$ as a renormalized average over all these chambers.
Note that in this case $X_S^F(C,d)$ need not be a Laurent polynomial
in  $y$ and  $E^F_S(C,d)$, and
$S_S^F(C,d)$ need not be integers. 
\end{enumerate}
\end{rem}

\begin{rem} Note that we did not define 
$\xx^{S,F}_C$ for $F\in \SS_S$ in case $K_SF=0$ and $FC$ even. 
The point is that in this
case
$\Theta^{H,F}_{\Gamma,C}(2\tau,K_Sz)$ is not well-defined:
$\Theta^{H,F}_{\Gamma,C}(2\tau,\cdot)$ has a pole along $F^\perp$.
If $S$ is the blowup of $\P_2$ in $9$ points and $F=-K_S$,
it is easy to see that $H(M^L_S(C,d))$ is constant for $L$ near $F$, and one can
therefore define  $X_y(M^F_S(C,d)):=X_y(M^L_S(C,d))$. This case has been studied by
Yoshioka \cite{Y4} in order to check predictions from \cite{M-N-V-W}.

In future, whenever we deal with $\xx^{S,F}_C$, $\ee^{S,F}_C$ for 
$F\in\SS_S$, we implicitely assume that  $K_SF\ne 0$ or $CF$ is odd. 
\end{rem}

\begin{cor}\label{crosscor} Let $F$, $G$ in $\overline \cc_S$, then 
\begin{align*}
\xx^{S,F}_C-\xx^{S,G}_C&=(y^{\frac{1}{2}}-y^{-\frac{1}{2}})
\frac{\eta(\tau)^{2\sigma(S)-2}}
{\theta_{1,1}(\tau,z)^2}
\Theta^{G,F}_{\Gamma,C}(2\tau,K_Sz),\\
(\ee^{S,F}_C)_0-(\ee^{S,G}_C)_0&=\text{\rm Coeff}_{2\pi i z}
\big(\Theta^{G,F}_{\Gamma,C}(2\tau, K_Sz)\big),\\
\ee^{S,F}_C-\ee^{S,G}_C
&=\frac{1}{\eta(\tau)^{2e(S)}}
\hbox{\rm Coeff}_{2\pi iz}\big(\Theta^{G,F}_{\Gamma,C}(2\tau, K_Sz)\big),\\
\Sigma^{S,F}_C-\Sigma^{S,G}_C&=
\frac{\eta(\tau)^{2\sigma(S)}}{2i\eta(2\tau)^{4}}
\Theta^{G,F}_{\Gamma,C}(2\tau,K_S/2).
\end{align*}
\end{cor}
\begin{pf} This is straightforward from Thm.~\ref{walltheta},
Def.~\ref{extend} and the cocycle condition (\ref{cocycle}).
\end{pf}

\subsection{Birational properties}
Let $\widehat S$ be the blowup of $S$ in a point, and let $E$ be the 
exceptional divisor.
We identify $H^2(S,\Z)$ with $E^\perp\subset H^2(\widehat S,\Z)$.

\begin{cor}(Blowup formulas)\label{blowcor}
Assume $C\not\in 
2H^2(S,\Z)$, then for all $F\in \overline \cc_S$ we have
\begin{enumerate}
\item
$\displaystyle{
\xx^{\widehat S,F}_C=
\frac{\theta_{0,0}(2\tau,z)}{\eta(\tau)^2}\xx^{S,F}_C,\quad
\xx^{\widehat S,F}_{C+E}
=\frac{\theta_{1,0}(2\tau,z)}{\eta(\tau)^2}\xx^{S,F}_C,}$
\item 
$\displaystyle{\Sigma^{\widehat S,F}_C=\frac{1}{\eta(2\tau)}\Sigma^{S,F}_C,\qquad
\Sigma^{\widehat S,F}_{C+E}=0.}$
\item
If $F\in \cc_S$ or $F\in\delta(\cc_S)$  and $FC$ is odd,
then 
$$\ee^{\widehat S,F}_C=\frac{\theta(2\tau)}{\eta(\tau)^2}\ee^{S,F}_C,\quad
\ee^{\widehat S,F}_{C+E}=\frac{\theta_{1,0}^0(2\tau)}{\eta(\tau)^2}\ee^{S,F}_C.
$$
\end{enumerate}
\end{cor}

\begin{pf} If $C\not \in 2H^2(S,Z)$ and $F\in \cc_S^G$ does not lie 
on a wall of type
$(C)$, or if $F\in \delta(\cc^G_S)$, and $CF$ is odd,
 then by \cite{Be}
all the cohomology  of $M^F_S(C,d)$, $M^F_{\widehat S}(C,d)$ and of 
$M^F_{\widehat S}(C+E,d)$ is of type $(p,p)$. Therefore in this case the result follows 
from the blowup formulas from \cite{L-Q1},
\cite{L-Q2}. Alternatively  one can use Prop.~\ref{blowprop}.
 In order to check the result for general $F$, we have to 
check that the blowup formula is compatible with our extension.
Let $\Gamma:=H^2(S,\Z)$ and $\widehat \Gamma:=H^2(\widehat S,\Z)$ with the negative
of the intersection forms.  By definition the compatibility of $\xx^{S,F}_C$
with the blowup formulas amounts to   the easy formulas
\begin{align*}
\Theta^{G,F}_{\widehat
\Gamma,C}(\tau,K_{\widehat S}z)&=\theta_{0,0}(\tau,z)\Theta^{G,F}_{
\Gamma,C}(\tau,K_Sz),\\
\Theta^{G,F}_{\widehat
\Gamma,C+E}(\tau,K_{\widehat S}z)&=\theta_{1,0}(\tau,z)\Theta^{G,F}_{
\Gamma,C}(\tau,K_Sz).
\end{align*}
The result for the signatures follows by
$$\frac{\theta_{0,0}(2\tau,1)}{\eta(\tau)^2}=\frac{\theta_{0,1}^0(2\tau)}{\eta(\tau)^2}
=\frac{1}{\eta(2\tau)},\quad
\theta_{1,0}(2\tau,1)=\theta_{1,1}(\tau,0)=0.
$$
\end{pf}

Immediately from (2) we get:
\begin{cor}
Assume that $\pi:S\to X$ is the blowup of a rational surface $X$ in finitely many 
points.  Let $C\in H^2(S,\Z)\setminus \pi^*(H^2(X,\Z))$. Let $F\in
\pi^*(\overline \cc_X^G)$ not on  a wall of type $(C,d)$. If $F\in
\pi^*(\delta(\cc_X^G))$, assume that
$CF$ is odd. Then  $\sigma(M^F_S(C,d))=0$.
\end{cor}

Cor.~\ref{crosscor} expresses the differences 
$\xx^{S,F}_C-\xx^{S,G}_C$ in terms of the theta functions
$\Theta^{F,G}_{\Gamma,C}$. We now want to show that
 in case $C\not \in H^2(S,\Z)$ we can for suitable $G$ also express
 $\xx^{S,F}_C$, $\ee^{S,F}_C$ and $\Sigma^{S,F}_C$ in terms of 
$\Theta^{F,G}_{\Gamma,C}$.
We use the following easy fact (see e.g.\cite{Q1}, \cite{H-L}).
\begin{lem}\label{vanish}
Let $\pi:X\to \P_1 $ a rational ruled surface. Let $S$ be obtained from $X$
by successively blowing up a number of points.    Let $G$ be pullback of the class of a
fibre of $\pi$. Let $C\in H^2(S,\Z)$ with $CG$ odd. Then $M^G_S(C,d)=\emptyset$ for all
$d$.
\end{lem}

\begin{prop}\label{ratprop} Let $F\in \overline \cc_S$.
Assume $C\not \in 2H^2(S,\Z)$. There exists a blowup $\widetilde S$ of $S$ in $n$
points and 
a $G\in \SS_{\widetilde S}^G$ such that $\xx^{\widetilde S,G}_C=0$.
Let $E_1,\ldots,E_n$ be the classes of the exceptional divisors,
and write $\widetilde \Gamma:=\Gamma\oplus\<E_1,\ldots,E_n\>$.
Then 
\begin{align*}
\xx^{S,F}_C&=\frac{
\eta(\tau)^{2\sigma(S)-2}}{\theta_{1,1}(\tau,z)^2\theta_{0,0}(2\tau,z)^n}
\Theta^{G,F}_{\widetilde\Gamma,C}
(2\tau,K_{\widetilde S}z)
,\\
\ee^{S,F}_C&=\frac{\hbox{\rm Coeff}_{2\pi i z}\big(\Theta^{G,F}_{\widetilde
\Gamma ,C} (2\tau,K_{\widetilde S}z)\big)}
{\eta(\tau)^{2e(S)}\theta(2\tau)^n},\\
\Sigma^{S,F}_C&=\frac{\eta(\tau)^{2\sigma(S)}}{2i\eta(2\tau)^4}
\frac{\Theta^{G,F}_{\widetilde \Gamma,C}
(2\tau,K_{\widetilde S}/2)}{\theta_{0,1}^0(2\tau)^n}
\end{align*}
\end{prop}
\begin{pf}
Any rational  surface $S$ can be blown up in such a way that 
$\widetilde S$ is a blowup of a ruled surface $X$, and  $CG$ is odd for 
$G$ the pullback of the fibre. Then  by Lem.~\ref{vanish} 
$M^G_S(C,d)=\emptyset$ for all $d$, 
and therefore $\xx^{\widetilde S,G}_C=0$.
The formulas are then a straightforward application
of Thm.~\ref{walltheta} and the blowup formula Cor.~\ref{blowcor}.
\end{pf}

\begin{cor}
$\xx^{S,H}_C$ is invariant under deformations of the triple $(S,H,C)$.
\end{cor}

Let $p_1,p_2,p_3$ be three non-collinear points in
$\P_2$.  Let $L_1,L_2,L_3$ be the lines through pairs of the  $p_i$ with
$p_i,p_j\in L_k$ for distinct indices $i,j,k$. Let $X$ be the blowup
of $\P_2$ in $p_1,p_2,p_3$, let $E_1,E_2,E_3$ be the exceptional divisors
and $\overline E_1,\overline E_2,\overline E_3$ the strict transforms of $L_1,L_2,L_3$. 
They can be blown down to points $\overline p_1,\overline p_2,\overline p_3$ to
obtain another projective plane $\overline\P_2$. Let $H$ and $\overline H$
be the hyperplane classes on $\P_2$ and $\overline \P_2$.
Let $S$ be the blowup of $X$ in additional points $p_4,\ldots, p_r$
with exceptional divisors $E_4,\ldots,E_r$. We denote the pullbacks of $H$ and
$\overline H$ by the same letters. Then
$\overline H=2H-E_1-E_2-E_3$, $\overline E_i=H-E_j-E_j$ for $i,j,k$ distinct in
$\{1,2,3\}$. We can view $S$ both as a blowup of $\P_2$ in $p_1,\ldots,p_r$ and as a
blowup of $\overline \P_2$ in $\overline p_1,\overline p_2,\overline p_3,
p_4,\ldots,p_r$. The change of viewpoint  amounts to a Cremona transform on $H^2(S,\Z)$
sending
$dH-\sum_{i=1}^r a_i E_i$ to 
$$
(2d-a_1-a_2-a_3)H-(d-a_2-a_3)E_1-(d-a_1-a_3)E_2-(d-a_1-a_2)E_3-\sum_{i=4}^r a_i
E_i.$$
This shows:
\begin{cor}\label{cremcor} Let $C\in H^2(S,\Z)$ and $F\in\overline \cc_S$.
Let $\G$ be the group generated by the Cremona transforms and the permutations 
of $E_1,\ldots,E_r$. Then for all $g\in \G$ we have
$\xx^{S,F}_C=\xx^{S,g(F)}_{g(C)}$, and, if $F\in \cc_S$, or $F\in \delta( \cc^G_S)$
with $FC$ is odd, then
$M_S^F(C,d)=M_S^{g(F)}(g(C),d)$ for all $d$.
\end{cor}

\section{Transformation properties on the boundary}

Vafa and Witten \cite{V-W} made predictions for the modular behaviour of   generating
functions of the Euler numbers of moduli spaces of sheaves on algebraic surfaces. Up to
eventual quasimodularity  their generating function $Z_C$ (which can be essentially
identified with $\ee^{F,S}_C$) should fullfil the equations
\begin{align}\label{VaWi} Z_C^S(\tau+1)&=\epsilon Z_C^S;\qquad Z_C^S(-1/\tau)= \pm
2^{-b_2(S)/2}\left(\frac{\tau}{i}\right)^{-e(S)/2}
\sum_{D}(-1)^{DC} Z_D^S.
\end{align} Here $\epsilon$ is a  root of unity, and 
$D$ runs through a system of representatives of $H^2(S,\Z)$ modulo $2H^2(S,\Z)$.

We want to show that for $F\in\SS_S$ a similar transformation behaviour holds for
$\xx^{S,F}_C$. Formulas similar to those of (\ref{VaWi}) for the Euler numbers then
follow as a corollary. In addition we also get the modular behaviour for the signatures.
For the purpose of this section we will for $F\in \overline \cc_S$ define 
$$
\xx^{S,F}_0:=\frac{\eta(\tau)^2}{\theta_{1,0}(2\tau,z)}\xx^{\widehat S,F}_E,\qquad
\ee^{S,F}_0:=\frac{\eta(\tau)^2}{\theta^0_{1,0}(2\tau)}\ee^{\widehat S,F}_E,
$$ where $\widehat S$ is the  blowup of $S$ in a point, and $E$ is the class of the
exceptional divisor, i.e. we formally use the blowup formulas Cor.~\ref{blowcor} (which
do not apply). This is similar to the approach for the  Donaldson invariants. We put
\begin{align*}
Y^{S,F}_C(\tau,z)&:=\frac{\xx^{F,S}_C}{y^{1/2}-y^{-1/2}},\\
F^{S,F}_C(\tau)&:=\ee^{S,F}_C(\tau)-\frac{K_S^2G_2(\tau)}{2\eta(\tau)^{2e(S)}}
\hbox{\rm Coeff}_{(2\pi i z)^{-1}}\Theta^{G,F}_{\Gamma,C}(2\tau,K_Sz).
\end{align*} Note that in case $CF$ odd we just have
$F^{S,F}_C(\tau):=\ee^{S,F}_C(\tau)$, because $\Theta^{G,F}_{\Gamma,C}(2\tau,K_Sz)$ is
holomorphic at $z=0$.

\begin{thm} For all $C\in H^2(S,\Z)$ and all $F\in \SS_S$ with $FK_S\ne 0$  we have
\begin{enumerate}
\item $Y^{S,F}_C$ transforms according to the rules
      \begin{align*} Y^{S,F}_C(\tau+1,z)&=(-1)^{-e(S)/6-C^2/2}
Y^{S,F}_C(\tau,z),\\ 
Y^{S,F}_C(-1/\tau,z/\tau)&=-i \sqrt{\frac{2\tau}{i}}^{-b_2(S)}
\exp(-\pi i (K_S^2/2+2)z^2/\tau)\\
&\qquad\cdot \left(\sum_{D\in H^2(S,\Z)/2H^2(S,\Z)}
(-1)^{CD}Y^{S,F}_D(\tau,z)\right).
\end{align*}
\item $F^{S,F}_C$ transforms according to 
\begin{align*} F^{S,F}_C(\tau+1)&=(-1)^{-e(S)/6-C^2/2}F^{S,F}_C(\tau)\\
F^{S,F}_C(-1/\tau)&=-\sqrt{\frac{\tau}{i}}^{-e(S)}2^{-b_2(S)/2}
\sum_{D\in H^2(S,\Z)/2H^2(S,\Z)}(-1)^{CD}F^{S,F}_D(\tau).
\end{align*}
\item Assume $C^2$ is odd. Then 
$\displaystyle{\frac{\eta(\tau)^2}{\eta(2\tau)^{\sigma(S)}}\Sigma^{S,F}_C}$ is a modular
function on $\Gamma(2)$.
\end{enumerate}
\end{thm}

\begin{rem} By the fact that $G_2(\tau)+1/(8\pi \Im(\tau))$ transforms like a modular
form of weight $2$, we could also define 
$$F^{S,F}_C(\tau):=\ee^{S,F}_C(\tau)+\frac{K_S^2}{16\pi\Im(\tau)\eta(\tau)^{2e(S)}}
\hbox{\rm Coeff}_{(2\pi iz)^{-1}}\Theta^{G,F}_{\Gamma,C}(2\tau,K_Sz),$$
and get the same transformation behaviour in 2.
\end{rem}

\begin{pf}  (1) Let $F,G\in \SS_S$. We first want to show that (1) holds if we replace
$Y^{S,F}_C$ and $Y^{S,F}_D$  by 
$Y^{S,F}_C-Y^{S,G}_C$ and $Y^{S,F}_D-Y^{S,G}_D$. By Cor.~\ref{crosscor} we have
$$Y^{S,F}_C-Y^{S,G}_C=\frac{\eta(\tau)^{2\sigma(S)-2}}{\theta_{1,1}(\tau,z)^2}
\Theta^{G,F}_{\Gamma,C}(2\tau,K_Sz).$$  
By (\ref{etatrans}) and (\ref{theta11trans}), we
know that 
\begin{align*}\frac{\eta(\tau+1)^{2\sigma(S)-2}}{\theta_{1,1}(\tau+1,z)^2}&=
(-1)^{-e(S)/6}\frac{\eta(\tau)^{2\sigma(S)-2}}{\theta_{1,1}(\tau,z)^2},\\
\frac{\eta(-1/\tau)^{2\sigma(S)-2}}{\theta_{1,1}(-1/\tau,z/\tau)^2}&=
-\Big(\frac{\tau}{i}\Big)^{-b_2(S)}e^{-2\pi i z^2/\tau}
\frac{\eta(\tau)^{2\sigma(S)-2}}{\theta_{1,1}(\tau,z)^2}.
\end{align*}
Furthermore, putting $T(\tau,z):=\Theta^{G,F}_{\Gamma,C}(2\tau,K_Sz)$ we get by
(\ref{thetatrans}) 
\begin{align*}
T(\tau+1,z)&=(-1)^{-C^2/2}T(\tau,z),\\
T(-1/\tau,z/\tau)&= 
\Theta^{G,F}_{\Gamma,C}(-2/\tau,K_Sz/\tau)\\
&= i\sqrt{\frac{\tau}{2i}}^{b_2(S)}\exp(-\pi i
K_S^2z^2/{2\tau}) 
\Theta^{G,F}_{\Gamma}(\tau/2,K_Sz/2+C/2).
\end{align*} 
Putting this together,
we obtain
$$
(Y^{S,F}_C-Y^{S,G}_C)(-1/\tau,z/\tau)=-i\sqrt{\frac{2\tau}{i}}^{-b_2(S)}\exp(-\pi i
(K_S^2/2+2)z^2/\tau)\Theta^{G,F}_{\Gamma}(\tau/2,K_Sz/2+C/2).
$$
Finally we have 
\begin{align}
\Theta^{G,F}_{\Gamma}(\tau/2,K_Sz/2+C/2)&=
\sum_{\xi\in
\Gamma}(\mu(\<\xi,G\>)-\mu(\<\xi,F\>)) q^{\<\xi,\xi\>/4}(-1)^{C\xi}y^{\<K_S,\xi\>/2}
\nonumber\\
&=\sum_{D}(-1)^{DC} \sum_{\xi\in
\Gamma+D/2}(\mu(\<\xi,G\>)-\mu(\<\xi,F\>)) q^{\<\xi,\xi\>}y^{\<K_S,\xi\>}\nonumber\\
&=\sum_{D}(-1)^{DC} \Theta^{G,F}_{\Gamma,D}(2\tau,K_Sz)\label{halflat}.
\end{align} This shows (1) for $Y^{S,F}_C-Y^{S,G}_C$. 
It is
therefore enough to show (1) $Y_C^{S,F}$ for all $S,C$ and   
 one particular  $F\in \SS_S$. 
 Let $\epsilon: \widehat
S\to S$ be the blowup in a point with exceptional divisor 
$E$. 
By the blowup formulas Cor.~\ref{blowcor} and the definition of
$\xx^{S,F}_0$ we get
$$Y^{\widehat S,F}_C(\tau,z)
=\frac{\theta_{0,0}(2\tau,z)}{\eta(\tau)^2}Y^{S,F}_C(\tau,z),
\quad 
Y^{\widehat S,F}_{C+E}(\tau,z)
=\frac{\theta_{1,0}(2\tau,z)}{\eta(\tau)^2}Y^{S,F}_C(\tau,z),$$
and therefore
\begin{align*}
Y^{\widehat S,F}_C(\tau+1,z)&
=(-1)^{-1/6}\frac{\theta_{0,0}(2\tau,z)}{\eta(\tau)^2}Y^{S,F}_C(\tau+1,z),\\
Y^{\widehat S,F}_{C+E}(\tau+1,z)&
=(-1)^{1/2-1/6}\frac{\theta_{1,0}(2\tau,z)}{\eta(\tau)^2}Y^{S,F}_C(\tau+1,z).
\end{align*}
By the transformation behaviour of $\theta_{0,0}$, $\theta_{1,0}$ (see
\cite{C} Sect. V.8.) and  $\eta$  we also see
\begin{align*}
Y^{\widehat S,F}_C(-1/\tau,z/\tau)&
=\sqrt{\frac{i}{2\tau}}\exp((\pi i/2)z^2/\tau)
\frac{\theta_{0,0}(\tau/2,z/2)}{\eta(\tau)^2}Y^{S,F}_C(-1/\tau,z/\tau),\\
Y^{\widehat S,F}_{C+E}(-1/\tau,z/\tau)&
=\sqrt{\frac{i}{2\tau}}\exp((\pi i/2)z^2/\tau)
\frac{\theta_{0,1}(\tau/2,z/2)}{\eta(\tau)^2}Y^{S,F}_C(-1/\tau,z/\tau).
\end{align*}
Using the   elementary identities
\begin{align*}\theta_{0,0}&(\tau/2,z/2)=\theta_{0,0}(2\tau,z)
+\theta_{1,0}(2\tau,z),\quad
\theta_{0,1}(\tau/2,z/2)=\theta_{0,0}(2\tau,z)
-\theta_{1,0}(2\tau,z)
\end{align*}
it follows that (1) holds for $Y^{S,F}_C$ for all $C\in H^2(S,\Z)$
if and only if it holds for $Y^{\widehat S,F}_{\overline C}$ for all 
$\overline C\in H^2(\widehat S,\Z)$.    As
any two rational surfaces can be connected by a sequence of  blowups and blow downs, it
is enough to check the result for $S=\P_1\times \P_1$  
and $F$ the class of a fibre of
the first projection. Let 
$G$ be the class of a fibre of the second projection. By Lem.~\ref{vanish} we have
$Y^{S,F}_G=Y^{S,F}_{F+G}=0$.  Denote by $\widehat \P_2$ the blowup of $\P_2$ in a point
with exceptional divisor $E_1$. Let $H\in H^2(\widehat \P_2,\Z)$ be the class of a
hyperplane. Let $\sigma:\widetilde \P_2\to \widehat \P_2$ be the blowup in a point with
exceptional divisor $E_2$. There exists a blowup 
$\epsilon:\widetilde \P_2\to \P_1\times\P_1$ with exceptional divisor $E$ such that
$\sigma^*(H-E_1)=\epsilon^*(F)$ and  $E=\epsilon^*(F)-E_2$. By Cor.~\ref{blowcor} and
our definition of $Y^{S,F}_0$ we get that 
\begin{align*} Y^{S,F}_F&=\frac{\eta(\tau)^2}{\theta_{0,0}(2\tau,z)}Y^{\widetilde
\P_2,\epsilon^*(F)}_{\epsilon^*(F)}= Y^{\widehat \P_2,H-E_1}_{H-E_1},\\
Y^{S,F}_0&=\frac{\eta(\tau)^2}{\theta_{1,0}(2\tau,z)}Y^{\widetilde \P_2,\eps^*(F)}_E=
Y^{\widehat \P_2,H-E_1}_{H-E_1}.
\end{align*} 
By $F^2=G^2=0$, $FG=1$,  (1)  follows.

(2) By an argument that is very similar to that at the end of the proof of (1),
 it is enough to show the 
formula for the difference $F^{S,F}_C-F^{S,G}_C$, for $F,G\in \SS_S$. 
The transformation
behaviour
$F^{S,F}_C(\tau+1)=(-1)^{-e(S)/6-C^2/2}F^{S,F}_C(\tau)$
follows immediately from the corresponding transformation behaviour of $Y^{S,F}_C$.
By the transformation behaviour of $\Theta^{G,F}_{\Gamma,C}(2\tau,K_Sz)$, 
the transformation properties (\ref{g2trans}) of $G_2$ and
(\ref{halflat}), we see that
$U_C(\tau,z):=\exp(2\pi^2K_S^2G_2(\tau)z^2)\Theta^{G,F}_{\Gamma,C}(2\tau,K_Sz)$,
transforms  according to 
$$U_C(-1/\tau,z/\tau)
=i\sqrt{\frac{\tau}{2i}}^{-b_2(S)}  \sum_{D\in
H^2(S,\Z)/2H^2(S,\Z)}U_D(\tau,z).$$
By definition 
$$(F^{S,F}_C-F^{S,G}_C)(\tau)=\frac{1}{\eta(\tau)^{2e(S)}}
\hbox{\rm Coeff}_{2\pi iz} U_C(\tau,z),$$
and (2) follows.

(3) Let $\widehat S$ be the blowup of $S$ in a point. By the blowup formula
Cor.~\ref{blowcor}, we see that the statements for $(S,F,C)$ and $(\widehat S,F,C)$ 
are equivalent. Therefore, by Prop.~\ref{ratprop}, we can assume that there exists
a $G\in\SS^G_S$ such that $\xx_C^{S,G}=0$ and 
$$\Sigma^{S,F}_C(\tau)=\frac{\eta(\tau)^{2\sigma(S)}}{2i\eta(2\tau)^4}
\Theta^{G,F}_{\Gamma,C}(2\tau,K_S/2).$$
By \cite{G-Z}, Thm 3.13.1) the function
$$\tau\mapsto
G(\tau)=\frac{\theta(\tau)^{\sigma(S)}}{f(\tau)}\Theta^{G,F}_{\Gamma,C}(\tau,C/2)$$ is a
modular function on $\Gamma_u$. Therefore, by Rem.~\ref{gammau},
$\tau\mapsto G(2\tau+1)$ is a
modular  function on $\Gamma(2)$. By (\ref{thetatrans}), we have 
$$\Theta^{G,F}_{\Gamma,C}(2\tau+1,C/2)=
(-1)^{-3C^2/4+CK_S/2}\Theta^{G,F}_{\Gamma,C}(2\tau,K_S/2),$$
and by Rem.~\ref{gammau} we see that 
$$\frac{\theta(2\tau+1)^{\sigma(S)}}{f(2\tau+1)}=\frac{\eta(\tau)^{2\sigma(S)+2}}
{\eta(2\tau)^{\sigma(S)+4}}.$$
The result follows.\end{pf}

\section{The signature and the Donaldson invariants}

Let again $S$ be a rational algebraic surface, let $H\in \cc_S$, and let
$E$ be a differentiable complex vector bundle on $S$, with Chern classes 
$(C,c_2)$. Let $d:=c_2-C^2/4$ and  $e:=4d-3$.
Let $A_e(S)$ be the set of polynomials of weight $e$ in $H_2(S,\Q)\oplus H_0(S,\Q)$,
where $a\in H_2(S,\C)$ has weight $1$, and the class $p\in H_0(S,\Z)$ of a point has
weight $2$.
The Donaldson invariants corresponding to $E$, the Fubini-Study metric associated
to 
$H$, and the homology orientation determined by the connected component of 
$\big\{L\in H^2(S,\R)\bigm| L^2>0\big\}$ containing $H$ are a linear map
$\Phi^{S,H}_{C,e}:A_e(S)\to \Q$. 
Let 
$$\Phi^{S,H}_{C}:=\sum_{e\ge 0}\Phi^{S,H}_{C,e}:
A_*(S):=\bigoplus_{e\ge 0} A_e(S)\to \Q.$$
In \cite{K-M} it is shown (more generally for simply connected $4$-manifolds $S$ with
$b_+=1$, where now $H$ is the period point of a Riemannian metric on $S$), that
$\Phi^{S,H}_{C,e}$ depends only on the chamber of type
$(C,d)$ of 
$H$, and that $\Phi^{S,H}_{C,e}-\Phi^{S,L}_{C,e}$ can be expressed as a sum
of wallcrossing terms 
$\delta^S_{\xi,e}$, for $\xi$ running through the classes of type $(C,d)$
with $\xi H<0<\xi L$. Kotschick and Morgan
make a conjecture about the structure of the $\delta^S_{\xi,e}$.

\begin{conj}  \cite{K-M}\label{KMconj}
$\delta^S_{\xi,e}(x^e)$ is a polynomial in
$\xi x$, $x^2$ whose coefficients depend only on $\xi^2$, $e$ and the 
homotopy type of $S$.
\end{conj}

Using this conjecture the difference $\Phi^{S,H}_{C,e}-\Phi^{S,L}_{C,e}$
was in \cite{Go3} and \cite{G-Z} expressed in terms of modular forms and 
theta functions.
In a series of (in part forthcoming) papers \cite{F-L1}, \cite{F-L2}, 
\cite{F-L3}, \cite{F-L4}
Feehan and Leness work towards a proof of conjecture \ref{KMconj}.
In \cite{F-L1} some necessary gluing results are proven.

We will in this section assume Conjecture \ref{KMconj}.
and show that  for any class $H\in \overline \cc_S$ the generating function for
the signatures
$\sigma(M_S^H(C,d))$ is also (with respect to a different development parameter)
the generating function for the Donaldson invariants $\Phi^{S,H}_{C}(p^r)$,
evaluated on the powers of the point class $p$.
The reason for this result is that both Donaldson invariants and the signatures
of the moduli spaces vanish in certain chambers, the
chamber structures for  Donaldson invariants and signatures are the same, and the 
wallcrossing terms for Donaldson invariants and signatures are related.

\begin{thm}\label{donsig} Assume Conjecture \ref{KMconj}.
Let $H\in \overline \cc_S$. Then
$$\Phi^{S,H}_{C}(p^r)=(-1)^{(CK_S+1)/2}\hbox{\rm Coeff}_{\overline u(\tau)^{r+1}}
\Big[\frac{4\eta(\tau)^2}{\eta(2\tau)^{\sigma(S)}}\Sigma^{S,H}_C\Big].$$
Here  $\hbox{\rm Coeff}_{\overline u(\tau)^{r+1}}V(\tau)$ is the coefficient
of $\overline u(\tau)^{r+1}$ in the Laurent development of $V(\tau)$ in powers of 
$\overline u(\tau)$.
In particular, if $H\in \cc^G_S$ does not lie  on a wall of type $C$ or $F\in
\delta(\cc_S)$ with
$CF$ odd, then
$$\sum_{d\ge 0}(-1)^{e(d)/2}\sigma(M_S^H(C,d))q^{d-e(S)/12}=
(-1)^{(CK_S+1)/2}\frac{\eta(2\tau)^{\sigma(S)}}{4\eta(\tau)^2}
\left(\sum_{r\ge 0} \Phi^{S,H}_{C}(p^r) \overline u(\tau)^{r+1}\right).$$
\end{thm}
\begin{pf}
We note that the result is trivially true if $C^2$ is even. We assume
that $C^2$ is odd.

{\it Case 1:} Assume that $S$ is the blowup of a ruled surface and that $CG$ is 
odd for $G$ the pullback of the class of the fibre of the ruling.
Then we get by Prop.~\ref{ratprop}
$$\Sigma^{S,H}_C=\frac{\eta(\tau)^{2\sigma(S)}}{2i\eta(\tau)^4}
\Theta^{G,H}_{\Gamma,C}(2\tau,K_S/2).$$
On the other hand we get by \cite{G-Z} Cor.~4.3 and Lem.~5.1
$$\Phi^{X,H}_C(p^r)=\hbox{\rm Coeff}_{u(\tau)^{r+1}}
\left[(-1)^{\frac{3}{4}C^2}\frac{2\theta(\tau)^{\sigma(S)}}{f(\tau)}
\Theta^{G,H}_{\Gamma,C}
(\tau,C/2)\right].$$
We make the transformation $\tau\to 2\tau+1$.
By (\ref{thetatrans}) we get 
$$\Theta^{G,H}_{\Gamma,C}(2\tau+1,C/2)
=(-1)^{-3C^2/4+CK_S/2}\Theta^{G,H}_{\Gamma,C}(2\tau,K_S/2).$$
Using also Rem.~\ref{gammau}, we get 
\begin{align*}
\Phi^{S,H}_C(p^r)&=(-1)^{CK_S/2}
\hbox{\rm Coeff}_{\overline u(\tau)^{r+1}}\left[
\frac{2\eta(\tau)^{2\sigma(S)+2}}{\eta(2\tau)^{\sigma(S)+4}}
\Theta^{H,G}_{\Gamma,C}(2\tau,K_S/2)\right]\\
&=(-1)^{(CK_S+1)/2}\hbox{\rm Coeff}_{\overline
u(\tau)^{r+1}}\Big[\frac{4\eta(\tau)^2}{\eta(2\tau)^{\sigma(S)}}\Sigma^{S,H}_C\Big].
\end{align*}
This shows the first part. To show the second part, we need to see that the smallest
power of $\overline u(\tau)$ that occurs in the development of 
$\frac{4\eta(\tau)^2}{\eta(2\tau)^{\sigma(S)}}\Sigma^{S,H}_C$ is $\overline u(\tau)$. 
We see that $\overline u(\tau)$ is $q^{\frac{1}{2}}$ multiplied with a power
series in $q$. In case $C^2\equiv 1$ modulo $4$ it follows from the definition that 
$\Sigma^{S,H}_C$ is
$q^{-e(S)/12+3/4}$ multiplied with a power series in $q$. If $C^2\equiv 3$ modulo $4$,
the fact that $M^S_H(C,d)$ is only nonempty if the expected dimension  $4d-3$ is
nonnegative implies that  is $q^{-e(S)/12+5/4}$ multiplied with a power
series in $q$. This shows the second part.

{\it General case:}
Let $\widetilde S$ be the blowup of $S$ in $n$ points, so that
case 1 applies to $\widetilde S$.
Then by the blowup formulas for the Donaldson invariants \cite{F-S} and by
Cor.~\ref{blowcor}
$$\Phi^{S,H}_C(p^r)=\Phi^{\widetilde S,H}_{C}(p^r), \quad
\Sigma^{\widetilde S,H}_C=\frac{1}{\eta(2\tau)^n}\Sigma^{S,H}_C.$$
The result follows.
\end{pf}

\begin{cor} If $F\in \SS^G_S$, $CF$ is odd and $\sigma(S)>-8$, then
$\sigma(M_S^F(C,d))=0$.
\end{cor}
\begin{pf}
This follows immediately from \cite{G-Z} Cor.~5.5.
\end{pf}

\begin{rem}
Thm.~\ref{donsig} and the results for the $K3$ surface suggest that
there should be a general formula relating the Donalson invariants and
the signatures of the moduli spaces $M^H_S(C,d)$ for all simply connected 
algebraic surfaces
$S$ even if $p_g(S)>0$. In general the moduli spaces 
$M^H_S(C,d)$ will be very singular, and one first has to find a suitable definition
of the signature.
The simplest formula that fits the known data seems to be the following:
\begin{align*}\sum_d &(-1)^{e(d)/2} \sigma(M^H_S(C,d))q^{d-e(S)/12}=\pm
\frac{\eta(2\tau)^{\sigma(S)}}{(2\eta(\tau))^{2\chi(\oo_S)}}
w(\tau)\left(
\sum_{r\ge 0} \Phi^S_C(p^r)\overline u(\tau)^{r+1}\right)^{\chi(\oo_S)},
\end{align*}
where 
$$w(\tau)=\begin{cases}\frac{1}{2}\left(\left(\frac{2\overline
u(\tau)+1}{2\overline u(\tau)}\right)^{\chi(\oo_S)}+\left(\frac{2\overline
u(\tau)-1}{2\overline u(\tau)}\right)^{\chi(\oo_S)}\right), &\hbox{if
}3\chi(\oo_S)-C^2\equiv 2
\hbox{ mod }4,\\
\frac{1}{2}\left(1+2\overline u(\tau))^{\chi(\oo_S)}-(-1)^{\chi(\oo_S)}(1-2\overline 
u(\tau))^{\chi(\oo_S)}\right),&
\hbox{if }3\chi(\oo_S)-C^2\equiv 0
\hbox{ mod }4.\end{cases}
$$
The formula has the following features:
\begin{enumerate}
\item It gives the correct result for rational surfaces and for $K3$-surfaces.
\item It is compatible with the blow-up formulas of \cite{F-S} for the Donaldson
invariants and with those of Prop.~\ref{blowprop} for the signatures.
\item It is compatible with taking the disjoint union of algebraic surfaces.
\end{enumerate}
\end{rem}
The formula for the rational surfaces is just Thm.~\ref{donsig}, and the compatibility
with the blowup formulas is obvious. We check the formula for the $K3$ surfaces.
We know by \cite{G-H} that, for generic polarization $H$ and suitable $C\in H^2(X,\Z)$,
the moduli space
$M^H_X(C,d)$ has the same Hodge numbers as $X^{[2d-3]}$. Let $L$ and $M$ be two such
classes in $H^2(X,\Z)$, satisfying $L^2\equiv 2$ modulo $4$ and
$M^2\equiv 0$ modulo $4$.  Using
(\ref{hilbhodge}), this gives 
$$\sum_{d\ge 0}(-1)^{e(d)}(\sigma(M^H_X(L,d))+\sigma(M^H_X(M,d))q^{d-2}=
\sum_{d\ge 0}(-1)^n\sigma(S^{[n]})q^{(n-1)/2}=\frac{1}{\eta(\tau)^4\eta(\tau/2)^{16}}.$$
For the Donaldson invariants we have the following results:
$X$ fullfils the simple type condition and 
$\Phi^{S,H}_C=(-1)^{C^2/2}$ for all $C\in H^2(S,\Z)$ (see e.g. \cite{Kr-M}).
Therefore we get 
$\Phi^{X,H}_L(p^{2r})=-2^{2r}$ and $\Phi^{X,H}_M(p^{2r+1})=2^{2r+1}$, i.e. 
$$\sum_{r\ge 0}(\Phi^{X,H}_L(p^r)+\Phi^{X,M}_L(p^r))\overline u(\tau)^{r+1}
=-\frac{\overline u(\tau)}{1+2\overline
u(\tau)}=4\frac{\eta(2\tau)^8}{\eta(\tau/2)^8}.$$ 
The last identity is an elementary
exercise in modular forms (e.g. one multiplies both sides with a suitable modular form
on $\Gamma(2)$ such that they both become modular forms on $\Gamma(2)$ and
compares the first few coefficients). Putting this together, we obtain
\begin{align*}
\sum_{d\ge 0}
(\sigma(M_X^H(L,d)-&\sigma(M_X^H(M,d))q^{d-2}=
\frac{\eta(2\tau)^{\sigma(X)}}{(4\eta(\tau))^{2\chi(\oo_X)}}
\left(\sum_{r\ge 0}(\Phi^{X,H}_L(p^r)+\Phi^{X,H}_M(p^r))\overline u(\tau)^{r+1}\right)^2
\\
&=\frac{\eta(2\tau)^{\sigma(X)}}{(4\eta(\tau))^{2\chi(\oo_X)}}(1-2\overline u(\tau))^2
\left(\sum_{r\ge 0}\Phi^{X,H}_L(p^r) \right)^2\\
&=\frac{\eta(2\tau)^{\sigma(X)}}{(4\eta(\tau))^{2\chi(\oo_X)}}(1-1/(2\overline
u(\tau)))^2
\left(\sum_{r\ge 0}\Phi^{X,H}_M(p^r) \right)^2.
\end{align*}
The result follows by collecting the odd powers of $\overline u(\tau)$ (for $L$)
and the even powers of $\overline u(\tau)$ (for $M$).

\begin{rem} Note that $\overline u(\tau)$ is the modular function  on
$\Gamma(2)$ that occurrs in a natural way in physics (\cite{W}, there it is
called $u$).  In \cite{M-W} the Donaldson invariants of $4$-manifolds with $b=1$
were (using physics arguments) related to (Borcherds type \cite{Bo}) integrals
over the "$u$-plane"
$\H/\Gamma(2)$. This suggests that also many results of this paper could be reformulated
in terms of such integrals.
\end{rem}

For the Euler number we can    prove  a weaker statement along the same 
lines. We can relate the generating functions for the difference of the Euler
numbers for two polarizations $H,L$ to the difference of certain Donaldson
invariants between $H$ and $L$. Let $k_S$ be the Poincar\'e dual of 
$K_S$.

\begin{prop} Let $H$, $L\in\cc_S$ not on a wall of type $0$. Then
\begin{align*}
&\ee^{S,H}_0-\ee^{S,L}_0=
\frac{i\eta(2\tau)^6}{2\theta(2\tau)^{\sigma(S)+2}\eta(\tau)^{e(S)}}
\left(\sum_{r\ge 0}\big(
\Phi^{S,H}_0(k_Sp^r)-\Phi^{S,L}_0(k_Sp^r)\big)u(2\tau)^{r+1}
\right).
\end{align*}
\end{prop}

\begin{pf} The proof is similar to the case of the signature. By Corollary
\ref{crosscor} we have 
$$
\ee^{S,H}_0-\ee^{S,L}_0=
\frac{1}{\eta(\tau)^{2e(S)}}\hbox{\rm Coeff}_{2\pi iz}
\left[\Theta^{L,H}_{\Gamma,0}(2\tau,K_S z)\right].$$ 
As $H,L\in\cc_S$, we see that 
$\hbox{\rm Coeff}_{2\pi iz}
\left[\Theta^{L,H}_{\Gamma}(2\tau,K_s z)\right]$ starts in degree $\ge 1/2$ in
$q$. Using this, we get from \cite{G-Z}, Cor 4.3
\begin{align*}
\sum_{r\ge 0}&
\big(\Phi^{S,H}_0(k_S p^r)-\Phi^{S,L}_0(k_S p^r)\big)u(2\tau)^{r+1}\\ 
&=\hbox{Coeff}_{2\pi i
z}\left[\frac{2\theta(2\tau)^{\sigma(S)}}{f(2\tau)^2}
\Theta^{H,L}_{\Gamma,0}(2\tau,K_S z)\right]\\
&=\frac{2i\theta(2\tau)^{\sigma(S)+2}\eta(\tau)^{2e(S)}}
{\eta(2\tau)^6}(\ee^{S,L}_0-\ee^{S,H}_0)
\end{align*}
\end{pf}

\begin{rem} This result can be reformulated as follows.
The expression
$$
\ee^{S,H}_0+
\frac{\eta(2\tau)^6}{2i\theta(2\tau)^{\sigma(S)+2}\eta(\tau)^{e(S)}}
\left(\sum_{r\ge 0}
\Phi^{S,H}_0(k_Sp^r)u(2\tau)^{r+1}
\right)
$$
is independent of  $H\in \cc_S$.
\end{rem}

\section{Examples}

\subsection{Rational ruled surfaces}

Let $S$ be a rational ruled surface. Let $F$ be the class of a fibre of
the ruling, and let $G$ be a section with $G^2\le 0$.
By Lem.~\ref{vanish} we know that 
$\xx^{S,F}_C=0$ if $CF=1$.
We will compute $\xx^{S,F}_F$ and $\ee^{S,F}_F$.
Furthermore we set
$M^{F_+}_S(F,d):=M^{F+\epsilon G}_S(F,d)$ for $\epsilon>0$ sufficiently small, so 
that there is no wall of type $(F,d)$ between $F$ and $F+\epsilon G$.

\begin{prop}
\begin{enumerate}
\item
$\displaystyle{
\xx^{S,F}_F=\frac{\big(y^{\frac{1}{2}}-y^{-\frac{1}{2}}\big)\eta(\tau)}
{\theta_{1,1}(\tau,z)^2\theta_{1,1}(\tau,2z)},}$

\item $\displaystyle{\sum_{d\ge 0}X_y(M^{F_+}_S(F,d))q^{d-\frac{1}{3}}=
\frac{y^{\frac{1}{2}}-y^{-\frac{1}{2}}}{\eta(\tau)^2\theta_{1,1}(\tau,z)^2}
\left(\frac{\eta(\tau)^3}
{\theta_{1,1}(\tau,2z)}
-\frac{1}{y-
y^{-1}}\right),}$

\item $\displaystyle{\ee^{S,F}_F=\frac{2G_2(\tau)}{\eta(\tau)^8}},\qquad$
$\displaystyle{\sum_{d\ge 0}e(M^{F_+}_S(F,d))q^{d-\frac{1}{3}}=
\frac{2G_2(\tau)+\frac{1}{12}}{\eta(\tau)^8}.}$
\end{enumerate}
\end{prop}

\begin{pf}
Let $F_1,F_2$ be the fibres of the two projections of
$\P_1\times\P_1$ to $\P_1$.
By a sequence of blowups and blowdowns  $(S,F)$ can be obtained from 
$(\P_1\times\P_1,F_1)$, where in each blowup $F$ is replaced by its total transform.
By the blowup formula Cor.~\ref{blowcor} we get
$\xx^{S,F}_F=\xx^{\P_1\times\P_1,F_1}_{F_1}$. We can therefore assume that
$S=\P_1\times\P_1$ and $F=F_1$, $G=F_2$. By Cor.~\ref{vanish} we get 
$\xx^{S,G}_F=0$ and
$$\xx^{S,F}_F=\frac{y^{\frac{1}{2}}-y^{\frac{1}{2}}}{\eta(\tau)^2\theta_{1,1}(\tau,
z)^2}\Theta^{G,F}_{\Gamma,F}(2\tau,-2Fz-2Gz).$$
By (\ref{thetaext}) we have
\begin{align*}\Theta^{G,F}_{\Gamma,F}(2\tau,x)&=
\eta(2\tau)^3\frac{\theta_{1,1}(\cdot,\<(F-G),\cdot\>)}{\theta_{1,1}(\cdot,-\<G,\cdot\>)
\theta_{1,1}(\cdot,\<F,\cdot\>)}
\Big|_{F/2}(2\tau,x)\\
&=\frac{\eta(2\tau)^3\theta_{0,1}(2\tau,\<(F-G),x\>)}{\theta_{0,1}(2\tau,-\<G,x\>)
\theta_{1,1}(2\tau,\<F,x\>)}.
\end{align*}
Thus 
$$
\Theta^{G,F}_{\Gamma,F}(2\tau,-2Fz-2Gz)=
\frac{\eta(2\tau)^3\theta^0_{0,1}(2\tau)}{\theta_{0,1}(2\tau,2z)\theta_{1,1}(2\tau,2z)}.
$$
By  (\ref{thetahalf}) we get 
$$\xx^{S,F}_F=\frac{\big(y^{\frac{1}{2}}-y^{-\frac{1}{2}}\big)\eta(\tau)}
{\theta_{1,1}(\tau,z)^2\theta_{1,1}(\tau,2z)}.$$
To get the $\chi_y$-genus of $M^{F^+}_S(F,d)$, we note that
$$\sum_{d\ge 0}X_y(M^{F^+}_S(F,d)q^{d-\frac{1}{3}}=
\frac{y^{\frac{1}{2}}-y^{-\frac{1}{2}}}
{\eta(\tau)^2\theta_{1,1}(\tau,z)^2}\big(\lim_{\epsilon\to 0}\Theta^{G,F+\epsilon
G}_{\Gamma,F}(\tau,-2Fz-2Gz)\big),$$
and by formula (3.9.1) from \cite{G-Z}
$$\Theta^{G,F}_{\Gamma,F}(\tau,-2Fz-2Gz)\big)-\lim_{\epsilon\to
0}\Theta^{G,F+\epsilon
G}_{\Gamma,F}(\tau,-2Fz-2Gz)\big)=\frac{1}{y-y^{-1}}.$$
To finally obtain the formulas for the Euler numbers 
we use the formula
$$\frac{\eta(\tau)^3}{\theta_{1,1}(\tau,z)}=\frac{1}{2\pi i z}\exp\left(
\frac{2}{k!} G_k(\tau)(2\pi i)^k\right),\quad (\hbox{see \cite{Z1}}),$$
and 
$\displaystyle{\hbox{\rm Coeff}_{2\pi iz}\frac{1}{y-y^{-1}}=-\frac{1}{12}.}$
\end{pf}

\subsection{The rational elliptic surface}

Let $m\in \Z_{\ge 0}$. Let $S$ be the blowup of $\P_2$ in $4m+5$ points,
and assume that
$F:=(m+2)H-mE_1-\sum_{i=2}^{4m+5}E_i$ is nef, e.g. 
$F$ is the fibre of a fibration of $S$ over $\P_1$, such that the genus
of the generic fibre is $m$.

\begin{thm}\label{2lpl1} If $m$ is odd, then
\begin{enumerate}
\item
$\xx^{S,F}_H+\xx^{S,F}_{E_1}=\displaystyle{
\frac{\theta_{1,1}(\tau,mz)}{\theta_{1,1}(\tau,z)\widetilde
\theta_{1,1}(\tau/2,z)\theta_{0,1}(\tau,(m-1)z)\eta(\tau)\eta(\tau/2)^{4m+3}},}$
\item $\displaystyle{\ee^{S,F}_H+\ee^{S,F}_{E_1} =\frac{m}{\eta(\tau/2)^{4m+8}},}$
\item $\Sigma^{S,F}_H+\Sigma^{S,F}_{E_1}
=\frac{1}{\eta(\tau)^2\eta(\tau/2)^{4m+4}}.$
\end{enumerate}
If $m$ is even, then 
\begin{enumerate}
\item
$\xx^{S,F}_{H+E_2}+\xx^{S,F}_{E_1+E_2}=\displaystyle{
\frac{\theta_{1,1}(\tau,mz)}{\theta_{1,1}(\tau,z)\widetilde
\theta_{1,1}(\tau/2,z)\theta_{0,1}(\tau,(m-1)z)\eta(\tau)\eta(\tau/2)^{4m+3}},}$
\item $\displaystyle{\ee^{S,F}_{H+E_2}+\ee^{S,F}_{E_1+E_2}
=\frac{m}{\eta(\tau/2)^{4m+8}},}$
\item 
$\Sigma^{S,F}_{H+E_2}=\Sigma^{S,F}_{E_1+E_2}=0.$
\end{enumerate}
\end{thm}

\begin{pf} We mostly deal with the case $m=2l-1$ odd. The proof in the case $m$ even is 
analogous. Let
$\Gamma=H^2(S,\Z)$ with the negative of the intersection form. Let
$G:=H-E_1$. Then by Lem.~\ref{vanish} and Cor.~\ref{crosscor}
$$\xx^{S,F}_H+\xx^{S,F}_{E_1}=
\frac{y^{\frac{1}{2}}-y^{-\frac{1}{2}}}{\eta(\tau)^{16l+2}\theta_{1,1}(\tau,z)^2}
\Big(\Theta^{G,F}_{\Gamma,H}(2\tau,K_Sz)+\Theta^{G,E_1}_{\Gamma,H}(2\tau,K_Sz)\Big).$$ 
Let $[G,F]$ be the lattice generated by $G $ and $F$, and let $[G,F]^\perp$ be its
orthogonal complement in $\Gamma$.  We write $\Lambda:=[G,F]\oplus [G,F]^\perp$. By
$\<F,F\>=\<G,G\>=0$, $\<F,G\>=-2$, $\<E_1,F\>=1-2l$, $ \<E_1,G\>=-1$, 
$\<E_2,F\>=-1$, $\<E_2,G\>=0$ we see that
$\Lambda$ has index $4$ in $\Gamma$, and that $0$, $E_1$, $E_2$, $E_2+E_1$ form a system
of representatives of 
$\Gamma$ modulo $\Lambda$.  Therefore we get by  (\ref{thetaext}):
\begin{align*}
\Theta^{G,F}_{\Gamma,H}&+\Theta^{G,F}_{\Gamma,E_1}=
\Theta^{G,F}_{\Lambda}\big(|_{0}+|_{E_1}+|_{E_2}+|_{E_1+E_2}\big)
\big(|_{H/2}+_{E_1/2}\big)\\
&=\Theta^{G,F}_{\Lambda}\big(|_{0}+|_{E_1}+|_{E_2}+|_{E_1+E_2}\big)
\big(|_{0}+|_{G/2}\big)|_{E_1}.
\end{align*}  Let $D_{8l}:=\big\{ (a_1,\ldots,a_{8l})\in \Z^{8l}\bigm| \sum_{i=1}^{8l}
a_i\hbox{ even}\big\}$.  Then the map
$\varphi:D_{8l}\to [G,F]^\perp$ defined by $(a_1,\ldots,a_{8l})\mapsto
\sum_{i=1}^{8l} a_{i}(E_{i+1}-G/2)$ is easily seen to be an isomorphism of lattices. It
is well-known (and easy to check) that 
$$\Theta_{D_{8l}}(\tau,(x_1,\ldots,x_{8l}))=\frac{1}{2}
\Big(\prod_{i=1}^{8l}\theta_{0,0}(\tau,x_i)
+\prod_{i=1}^{8l}\theta_{0,1}(\tau,x_i)\Big).$$  So we get by (\ref{thetaext1})
\begin{align*}\Theta^{G,F}_{\Lambda}(\tau,x)&=
\frac{\eta(2\tau)^3\theta_{1,1}(2\tau,\<F-G,x\>)}
{\theta_{1,1}(2\tau,\<F,x\>)\theta_{1,1}(2\tau,\<-G,x\>)}\cdot\\ &
\qquad\cdot\frac{1}{2}\Big(\prod_{i=2}^{8l+1}\theta_{0,0}(\tau,\<E_i-G/2,x\>)
+\prod_{i=2}^{8l+1}\theta_{0,1}(\tau,\<E_i-G/2,x\>)\Big).
\end{align*} If $H(\tau,x)$ is a  function $\H\times(\Gamma_\C)\to \C$  satisfying
$H(\tau,x)=\theta_{a,b}(n\tau,\<L,x\>)H_1(\tau,x)$ for some $L\in \Gamma_\Q$, then, for 
$W\in \Gamma$,
$$H|_W(\tau,x)=\theta_{a+2\<W,L/n\>,b}(n\tau,\<L,x\>)H_1|_W(\tau,x).$$  We have
$\<H,F\>=-(2l+1)$, $\<H,G\>=-1$ and $\<H,E_i\>=0$ for $i\ge 2$;
$\<E_1,F\>=-(2l-1)$, $\<E_1,G\>=-1$ and
$\<E_1,E_i\>=0$ for $i\ge 2$;
$\<E_2,F\>=-1$, $\<E_2,G\>=0$ and $\<E_2,E_2\>=1$, $\<E_2,E_i\>=0$ for $i\ge 3$.  We also
use repeatedly (\ref{plus2}). Using this we obtain the following: Put
\begin{align*} A(\tau,x)&:=\frac{\eta(2\tau)^3\theta_{1,1}(2\tau,\<F-G,x\>)}
{\theta_{1,1}(2\tau,\<F,x\>)\theta_{1,1}(2\tau,\<-G,x\>)},\
B(\tau,x):=\frac{-\eta(2\tau)^3\theta_{1,1}(2\tau,\<F-G,x\>)}
{\theta_{0,1}(2\tau,\<F,x\>)\theta_{0,1}(2\tau,\<-G,x\>)},\\
C(\tau,x)&:=\frac{\eta(2\tau)^3\theta_{0,1}(2\tau,\<F-G,x\>)}
{\theta_{1,1}(2\tau,\<F,x\>)\theta_{0,1}(2\tau,\<-G,x\>)},\
D(\tau,x):=\frac{\eta(2\tau)^3\theta_{0,1}(2\tau,\<F-G,x\>)}
{\theta_{0,1}(2\tau,\<F,x\>)\theta_{1,1}(2\tau,\<-G,x\>)},\\
\alpha(\tau,x)&:=\frac{1}{2}\Big(\prod_{i=2}^{8l+1}\theta_{0,0}(\tau,\<E_i-G/2,x\>)
+\prod_{i=2}^{8l+1}\theta_{0,1}(\tau,\<E_i-G/2,x\>)\Big),\\
\beta(\tau,x)&:=\frac{1}{2}\Big(\prod_{i=2}^{8l+1}\theta_{1,0}(\tau,\<E_i-G/2,x\>)
+\prod_{i=2}^{8l+1}\theta_{1,1}(\tau,\<E_i-G/2,x\>)\Big),\\
\gamma(\tau,x)&:=\frac{1}{2}\Big(\prod_{i=2}^{8l+1}\theta_{0,0}(\tau,\<E_i-G/2,x\>)
-\prod_{i=2}^{8l+1}\theta_{0,1}(\tau,\<E_i-G/2,x\>)\Big),\\
\delta(\tau,x)&:=\frac{1}{2}\Big(\prod_{i=2}^{8l+1}\theta_{1,0}(\tau,\<E_i-G/2,x\>)
-\prod_{i=2}^{8l+1}\theta_{1,1}(\tau,\<E_i-G/2,x\>)\Big).
\end{align*} Then we get 
\begin{align*}
\Theta^{G,F}_\Lambda|_{0}(\tau,x)&:=A(\tau,x)\alpha(\tau,x), 
\quad
\Theta^{G,F}_\Lambda|_{E_1}(\tau,x):=B(\tau,x)\beta(\tau,x),\\
\Theta^{G,F}_\Lambda|_{E_2}(\tau,x)&:=C(\tau,x)\gamma(\tau,x),
\quad
\Theta^{G,F}_\Lambda|_{E_1+E_2}(\tau,x):=D(\tau,x)\delta(\tau,x),\\
\Theta^{G,F}_\Lambda|_{G/2}(\tau,x)&:=D(\tau,x)\alpha(\tau,x), 
\quad
\Theta^{G,F}_\Lambda|_{G/2+E_1}(\tau,x):=C(\tau,x)\beta(\tau,x),\\
\Theta^{G,F}_\Lambda|_{G/2+E_2}(\tau,x)&:=B(\tau,x)\gamma(\tau,x),
\quad
\Theta^{G,F}_\Lambda|_{G/2+E_1+E_2}(\tau,x):=A(\tau,x)\delta(\tau,x).
\end{align*}  By $\<K_S,E_i-G/2\>=0$ and $\<E_1,E_i-G/2\>=1/2$ for $i\ge 2$, we see that
\begin{align*}
\alpha|_{E_1/2}(\tau,K_Sz)&=((\theta^0_{1/2,0}(\tau))^{8l}
+(\theta^0_{1/2,1}(\tau))^{8l})/2,\\
\beta|_{E_1/2}(\tau,K_Sz)&=((\theta^0_{3/2,0}(\tau))^{8l}
+(\theta^0_{3/2,1}(\tau))^{8l})/2,\\
\gamma|_{E_1/2}(\tau,K_Sz)&=((\theta^0_{1/2,0}(\tau))^{8l}
-(\theta^0_{1/2,1}(\tau))^{8l})/2,\\
\delta|_{E_1/2}(\tau,K_Sz)&=((\theta^0_{3/2,0}(\tau))^{8l}
-(\theta^0_{3/2,1}(\tau))^{8l})/2.
\end{align*} Now we note that by (\ref{plus2}) and (\ref{theta012})
$$\theta^0_{1/2,0}(\tau)=\theta^0_{3/2,0}(\tau) =\frac{\eta(\tau/2)^2}{\eta(\tau/4)},
\quad \theta^0_{1/2,1}(\tau)=-\theta^0_{3/2,1}(\tau).$$
 Putting things together, we get that 
$$(\Theta^{G,F}_{\Gamma,E_1}(\tau,K_Sz)+\Theta^{G,F}_{\Gamma,H}(\tau,K_Sz)=
(A+B+C+D)|_{E_1/2}(\tau,K_Sz) \cdot\frac{\eta(\tau/2)^{16l}}{\eta(\tau/4)^{8l}}.$$ The
orthogonal projections   of $0$,  $E_1$ and  $E_2$ and $E_1+E_2$ to $[F,G]$ are a system
of representatives of $[F/2,G/2]$ modulo $[F,G]$.   Therefore
\begin{align*}(A+B+C+D)(\tau,x) &=\Theta^{G,F}_{[G/2,F/2]}(\tau,x)\\
&=\frac{\eta(\tau/2)^3\theta_{1,1}(\tau/2,\<F/2-G/2,x\>)}
{\theta_{1,1}(\tau/2,\<F/2,x\>)\theta_{1,1}(\tau/2,\<-G/2,x\>)}.
\end{align*}  Finally by 
$\<E_1,F\>=1-2l$, $\<E_1,-G\>=1$, we get 
\begin{align*}
\Theta^{G,F}_{[G/2,F/2]}|_{E_1/2}(2\tau,K_Sz)&
=-\frac{\eta(\tau)^3\theta_{1,1}(\tau,\<F/2-G/2,K_S\>z)}
{\theta_{0,1}(\tau,\<F/2,K_S\>z)\theta_{0,1}(\tau,\<-G/2,K_S\>z)}\\
&=\frac{\eta(\tau)^3\theta_{1,1}(\tau,(2l-1) z)}
{\theta_{0,1}(\tau,(2-2l)z)\theta_{0,1}(\tau,z)}.
\end{align*}
 Putting this together, we obtain
\begin{align*}
\xx^{S,F}_H+\xx^{S,F}_{E_1}&=
\frac{(y^{\frac{1}{2}}-y^{-\frac{1}{2}})\eta(\tau)^3\theta_{1,1}(\tau,(2l-1)z)}
{\eta(\tau)^{16l+2}\theta_{1,1}(\tau,z)^2\theta_{0,1}(\tau,(2-2l)z)
\theta_{0,1}(\tau,z)}
\frac{\eta(\tau)^{16l}}{\eta(\tau/2)^{8l}}\\ &=
\frac{\theta_{1,1}(\tau,(2l-1)z)}{\theta_{1,1}(\tau,z)\widetilde
\theta_{1,1}(\tau/2,z)\theta_{0,1}(\tau,(2-2l)z)\eta(\tau)\eta(\tau/2)^{8l-1}}.
\end{align*}  In the last line we have used (\ref{thetahalf}). To get the formula for
$\ee^{S,F}_H+\ee^{S,F}_{E_1}$, take the limit $z\to 0$. It is immediate from
(\ref{theta11}) that
$\frac{\theta_{1,1}(\tau,(2l-1)z)}{\theta_{1,1}(\tau,z)}|_{z=0}=2l-1$ and
$\widetilde\theta_{1,1}(\tau/2,0)=\eta(\tau/2)^3$.  Therefore   
$$\ee^{S,F}_H+\ee^{S,F}_{E_1}
=\frac{2l-1}{\eta(\tau)\eta(\tau/2)^{8l+2}\theta_{0,1}^0(\tau)}
=\frac{2l-1}{\eta(\tau/2)^{8l+4}}.
$$ To get the formula for $\Sigma^{S,F}_H+\Sigma^{S,F}_{E_1}$ we put $z=\pi i$, to obtain
$$\Sigma^{S,F}_H+\Sigma^{S,F}_{E_1}=
\frac{1}{
\theta^0_{1,0}(\tau/2)\theta_{0,1}^0(\tau)\eta(\tau)\eta(\tau/2)^{8l-1}}=
\frac{1}{\eta(\tau)^2\eta(\tau/2)^{8l}}.
$$

In the case $m=2l$ even, we again have that
$0$, $E_1$, $E_2$, $E_1+E_2$ form a system of representatives of $\Gamma$ modulo
$\Lambda$. 
So we get
\begin{align*}
\Theta^{G,F}_{\Gamma,H+E_2}+\Theta^{G,F}_{\Gamma,E_1+E_2}
&=\Theta^{G,F}_{\Lambda}\big(|_{0}+|_{E_1}+|_{E_2}+|_{E_1+E_2}\big)
\big(|_{0}+|_{G/2}\big)|_{E_1+E_2}.
\end{align*}
Essentially the same computations as in the case $m$ odd give the result.
\end{pf}

Let now $S$ be the blowup of $\P_2$ in $9$ points. Let $H$ be the 
pullback of the hyperplane class, and let $E_1,\ldots,E_9$ the classes of the 
exceptional divisors. Let $F:=3H-\sum_{i=1}^9E_i$. Then $K_S=-F$. 
An interesting case is when  $S$
is a rational elliptic surface, and $F$ is the class of a fibre. 
In \cite{M-N-V-W} the generating functions of the Euler numbers $e(M_S(C,d))$ are 
predicted in case $CF$ is even.  This prediction was proven in \cite{Y4}.
As an immediate consequence of Thm.~\ref{2lpl1} we can compute the 
Hodge numbers of the $M^F_S(C,d)$ in case $FC$ is odd. For the Betti numbers this
result  was already obtained  (more generally for regular elliptic surfaces) in
\cite{Y6} using completely different methods. By \cite{Be} the result
about the Hodge numbers for $S$ is an direct consequence.

\begin{thm}\label{elliptic} Let $C\in H^2(S,d)$ with $C^2$ odd. Then $M_S^F(C,d)$ has
the same Hodge numbers as $S^{[2d-3/2]}$. In particular the Hodge numbers depend only on
$d$, and   we have
\begin{enumerate}
\item $\displaystyle{\sum_{d\ge 0}\big(X_y(M_S^F(H,d))+X_y(M_S^F(E_1,d))\big)q^{d-1}
=\frac{1}{\widetilde \theta_{1,1}(\tau/2,z)\eta(\tau/2)^9},}$
\item $\displaystyle{\sum_{d\ge 0}\big(e(M_S^F(H,d))+e(M_S^F(E_1,d))\big)q^{d-1}
=\frac{1}{\eta(\tau/2)^{12}},}$
\item $\displaystyle{\sum_{d\ge
0}\big(\sigma(M_S^F(H,d))-\sigma(M_S^F(E_1,d))\big)q^{d-1}
=\frac{1}{\eta(\tau)^2\eta(\tau/2)^{8}}.}$
\end{enumerate}
\end{thm}

\begin{rem} 
\begin{enumerate}
\item We can recover the Hodge numbers of $M_S^F(C,d)$:
\begin{align*}
\sum_{d\ge 0}X_y(M_S^F(H,d))q^{d-1}&=\frac{1}{\widetilde
\theta_{1,1}(\tau/2,z)\eta(\tau/2)^9}+\frac{i}{\widetilde
\theta_{1,1}((\tau+1)/2,z)\eta((\tau+1)/2)^9}\\
\sum_{d\ge 0}X_y(M_S^F(E_1,d))q^{d-1}&=\frac{1}{\widetilde
\theta_{1,1}(\tau/2,z)\eta(\tau/2)^9}-\frac{i}{\widetilde
\theta_{1,1}((\tau+1)/2,z)\eta((\tau+1)/2)^9}
\end{align*}

\item We can also use Thm.~\ref{donsig} to compute the generating functions
for the signatures. 
{}From \cite{G-Z} section 5.3 we get 
$\Phi^{S,F}_H(1+p/2)=1$, and by the simple type condition 
$\Phi^{S,F}_H(p^r)$ is $2^r$ if $r$ even and $0$ otherwise.
After some calculations this gives 
$$\sum_{d\ge 0}\sigma(M_S^F(H,d))q^{d-1}=12
e_2(2\tau)\frac{\eta(2\tau)^8}{\eta(\tau)^{22}}.$$
A similar calculation using \cite{G-Z} section 5.3 and Thm.~\ref{donsig}
gives 
$$\sum_{d\ge 0}\sigma(M_S^F(E_1,d))q^{d-1}=-8
\frac{\eta(2\tau)^{16}}{\eta(\tau)^{26}}.$$
It is an exercise in modular forms to show that 
$$
12e_2(2\tau)\frac{\eta(2\tau)^8}{\eta(\tau)^{22}}+8
\frac{\eta(2\tau)^{16}}{\eta(\tau)^{26}}=\frac{1}{\eta(\tau)^2\eta(\tau/2)^{8}},
$$
and thus to recover part (3) of Thm.~\ref{elliptic}.

\item If $X$ is a $K3$ surface, $L$ a primitive  line bundle  and $H$ a generic
ample line bundle on $X$, then it was shown in \cite{G-H} that 
$M^H_X(L,d)$ has the same Hodge numbers as $X^{[2d-3]}$, and in 
\cite{H} that $M^H_X(L,d)$ is deformation equivalent to $X^{[2d-3]}$.
There should be a similar proof of Thm.~\ref{elliptic} as that in \cite{G-H}.
Furthermore I expect that, in case $C^2$ odd, $M_S^F(C,d)$ is   deformation
equivalent to $S^{[2d-3/2]}$.
More generally similar results also should hold for arbitrary rank.

\item In physics the polarized rational elliptic surface $(S,F)$   is often called
$\frac{1}{2}K3$. This is related to the fact that one can degenerate an elliptic $K3$
surface to  the union of $2$ rational elliptic surfaces intersecting along  a fibre.
The generating function of the 
$\chi_y$-genera of the $M^H_X(L,d)$ ($L$ primitive and allowing $L^2$ both 
congruent $0$ modulo  $4$ and congruent $2$ modulo $4$) on the $K3$ surface is just the 
square of the generating function on $S$. One could ask whether this result can
also be shown by degenerating the moduli spaces $M^H_X(L,d)$. 
\end{enumerate}
\end{rem}

\begin{pf} (of Thm.~\ref{elliptic})
We first show  that the 
$M_S^F(C,d)$  depend only on $d$. Let $\G$ be the subgroup
of
$Aut(H^2(S,\Z))$ generated by the Cremona transforms and the permutations of 
$E_1,\ldots, E_{9}$. $F$ is invariant under the operation of $\G$,
and  therefore, by Cor.~\ref{cremcor},  $M_S^F(C,d)\simeq M_S^F(g(C),d)$ for all
$g\in\G$.   We can assume that
$C=nH-\sum_{i} a_i E_i$ with
$n,a_1,\ldots,a_{9}\in\{0,1\}$. Let $m$ be the number of indices $i\ge 1$ with
$a_i=1$.    By renumbering $E_1\ldots E_{9}$ we can assume that either
$C=H$ or $C=E_1$, in which case we are done, or 
$(h,a_1,a_2,a_3)$ is one of $(0,1,1,1)$ or   $(1,0,1,1)$.
The Cremona transform replaces
$(h,a_1,a_2,a_3)$ by 
$(1,0,0,0)$,   $(0,1,0,0)$,  and the result follows by induction on $m$.
As $K_SF\le 0$, the moduli spaces $M^F_S(C,d)$ are smooth, and by \cite{Be} all
their cohomology is of Hodge type $(p,p)$. Therefore  the
theorem follows from the  case
$l=1$ of Thm.~\ref{2lpl1}. \end{pf}

\end{document}